\documentclass[journal,twoside,web]
{IEEEtran}
\usepackage{cite}
\usepackage{amsmath,amssymb,amsfonts,amsthm}
\usepackage{graphicx}
\usepackage{algorithm}
\usepackage{hyperref}
\hypersetup{hidelinks=true}
\usepackage{textcomp}

\usepackage{bm}
\usepackage{algpseudocode} 
\usepackage{pifont}
\usepackage{caption}
\def\BibTeX{{\rm B\kern-.05em{\sc i\kern-.025em b}\kern-.08em
    T\kern-.1667em\lower.7ex\hbox{E}\kern-.125emX}}
\usepackage{atbegshi}%

 \usepackage{etex}
 
\usepackage{lipsum}
\usepackage{amsfonts}
\usepackage{graphicx}
\usepackage{epstopdf} 
\usepackage{enumerate} 
 \usepackage{amsmath}  
 \usepackage{subfig,bbm,enumerate,multirow,framed,blkarray}

\usepackage{verbatim,tcolorbox}
  \usepackage{tikz}
\usetikzlibrary{shapes,arrows}
\usetikzlibrary{positioning}  
\usepackage{pstricks-add} 
 \usepackage{wrapfig}
 \usepackage{pgfplots}
  
\usepackage{color}
\usepackage{savesym}
 
 \usepackage{graphicx}          

\algnewcommand{\algorithmicgoto}{\textbf{go to}}%
\algnewcommand{\Goto}[1]{\algorithmicgoto~\ref{#1}}%

\DeclareMathOperator*{\fsign}{fsign}

\usepackage{atbegshi}

\AtBeginDocument{\AtBeginShipoutNext{\AtBeginShipoutDiscard} \addtocounter{page}{-1}}

 \newenvironment{pf}{{\bfseries Proof.}}

\begin{document}

\newtheorem{theorem}{Theorem}[section] 
\newtheorem{corollary}[theorem]{Corollary}
\newtheorem{lemma}[theorem]{Lemma}
\newtheorem{proposition}[theorem]{Proposition}
\newtheorem{definition}[theorem]{Definition} 
\newtheorem{remark}[theorem]{Remark}
\newtheorem{assumption}[theorem]{Assumption} 
\newtheorem{example}[theorem]{Example}

\newcommand{\solnvi}{\text{SOL}(K,\mathbf{F})}
\newcommand{\solnncp}{{\text{NCP}}(\mathbf{F})}
\newcommand{\solnlcp}{{\text{LCP}}(\mathbf{y},\mathbf{M})}
\newcommand{\solnlelcp}{{\text{LeLCP}}(\mathbf{y},\mathbf{M})}

\newcommand{\anot}{\aaa_0}
\newcommand{\bnot}{\bbb_0}
\newcommand{\psib}{\pmb{\Psi}} 
\newcommand{\xib}{\pmb{\Xi}} 
\newcommand{\conv}{\text{conv}}
\newcommand{\plen}{\text{plen}}
 \newcommand{\fdir}{\mathbf{fdir}}
\newcommand{\qpstrong}{\mathcal{A}^+_{\text{QP}}}	
\newcommand{\qpweak}{\mathcal{A}^0_{\text{QP}}}
\newcommand{\qpinact}{\mathcal{A}^-_{\text{QP}}}
\newcommand{\qpact}{\mathcal{A}_{\text{QP}}}  
\newcommand{\cone}{P}  
\newcommand{\qp}{\text{QP}}
\newcommand{\lp}{\text{LP}}	 
   \newcommand{\strong}{\mathcal{A}^+}	
\newcommand{\weak}{\mathcal{A}^0}
\newcommand{\weakcrit}{\beta_{\text{crit}}}
\newcommand{\inact}{\mathcal{A}^-}
\newcommand{\act}{\mathcal{A}} 
\newcommand{\pset}{\mathcal{P}}
\newcommand{\iset}{\mathcal{I}}
\newcommand{\jset}{\mathcal{J}}
\newcommand{\nset}{\mathcal{N}}
\newcommand{\gset}{\mathcal{G}}
\newcommand{\sset}{\mathcal{S}}
\newcommand{\hset}{\mathcal{H}} 
	\algnewcommand\And{\textbf{and }}
	\algnewcommand\Or{\textbf{or }}
\newcommand{\pnot}{\pp_0}
\newcommand{\xnot}{\x_0}
\newcommand{\unot}{\uu_0}
\newcommand{\munot}{\mumu_0}
\newcommand{\etanot}{\etaeta_0} 
\newcommand{\rhonot}{\rhorho_0} 
\newcommand{\gamnot}{\gamgam_0} 
\newcommand{\gammanot}{\gamgam_0} 
\newcommand{\lambdanot}{\lambdalambda_0}
\newcommand{\dnot}{\dd_0} 
\newcommand{\ynot}{\mathbf{y}_0}
\newcommand{\znot}{\mathbf{z}_0} 
\newcommand{\vnot}{\mathbf{v}_0}

\newcommand{\minmin}{\pmb{\text{min}}}
\newcommand{\maxmax}{\pmb{\text{max}}}
\newcommand{\opt}{\pmb{\Phi}}
\newcommand{\mumu}{\pmb{\mu}}
\newcommand{\btheta}{\pmb{\theta}}
\newcommand{\nunu}{\pmb{\nu}}
\newcommand{\lambdalambda}{\pmb{\lambda}}
\newcommand{\deltadelta}{\pmb{\delta}}
\newcommand{\sx}{\widetilde{\pmb{x}}}
\newcommand{\smu}{\widetilde{\pmb{\mu}}}
\newcommand{\capmu}{\mathbf{U}}
\newcommand{\caplambda}{\mathbf{W}}
\newcommand{\lmin}{\mathbf{Lmin}}
\newcommand{\slmmin}{\mathbf{SLMmin}}
\newcommand{\lmmin}{\mathbf{LMmin}}
\newcommand{\lshiftmin}{\mathbf{SLmin}}
\newcommand{\lmax}{\mathbf{Lmax}}
\newcommand{\slmmax}{\mathbf{SLMmax}}
\newcommand{\lmmax}{\mathbf{LMmax}}
\newcommand{\lshiftmax}{\mathbf{SLmax}} 
\newcommand{\shift}{\mathbf{shift}}
\newcommand{\lmid}{\mathbf{Lmid}}
\newcommand{\lshiftmid}{\mathbf{SLmid}}
\newcommand{\coord}{\mathbf{e}} 
\newcommand{\E}{\mathbf{E}} 
\newcommand{\capw}{\mathbf{W}}
\newcommand{\caps}{\mathbf{S}}
\newcommand{\phimin}{\pmb{\Phi}_{\min}}
\newcommand{\phincp}{\pmb{\Phi}_{\text{NCP}}}
\newcommand{\psincp}{\pmb{\psi}_{\text{NCP}}}
\newcommand{\n}{\mathbf{n}}
\newcommand{\ximp}{\widetilde{\mathbf{x}}}
\newcommand{\zimp}{\widetilde{\mathbf{z}}}
\newcommand{\muimp}{\widetilde{\pmb{\mu}}}
\newcommand{\lambdaimp}{\widetilde{\pmb{\lambda}}}
\newcommand{\etaimp}{\widetilde{\pmb{\eta}}} 
\newcommand{\setmapsto}{\rvert\!\!\!\rightrightarrows}
\newcommand{\real}{\mathbb R}
\newcommand{\posint}{\mathbb{N}} 
\newcommand{\eps}{\epsilon}    
\newcommand{\tint}{[t_0,t_f]}
\newcommand{\fbar}{\hat{\mathbf{f}}_t}
\newcommand{\Fbar}{\hat{\mathbf{F}}_t} 
\newcommand{\gbar}{\hat{\mathbf{g}}_t}
\newcommand{\zx}{\pmb{\eta}_x}  
\newcommand{\zy}{\pmb{\eta}_y}
\newcommand{\zp}{\pmb{\eta}_p}
\newcommand{\etaxdot}{\pmb{\eta}_{\dot{x}}} 
\newcommand{\etaz}{\pmb{\eta}_z}  
\newcommand{\etazdot}{\pmb{\eta}_{\dot{z}}} 
\newcommand{\etaeta}{\pmb{\eta}}    
\newcommand{\rhorho}{\pmb{\rho}}  
\newcommand{\gamgam}{\pmb{\gamma}}
\newcommand{\zdot}{\mathbf{\dot{z}}}
\newcommand{\zz}{\mathbf{z}}
\newcommand{\w}{\mathbf{w}} 
\newcommand{\aaa}{\mathbf{a}}
\newcommand{\bbb}{\mathbf{b}}
\newcommand{\uu}{\mathbf{u}}
\newcommand{\vv}{\mathbf{v}}
\newcommand{\qq}{\mathbf{q}}
\newcommand{\zxdot}{\mathbf{z_{\dot{x}}}}   
\newcommand{\dx}{\mathbf{d_x}}   
\newcommand{\dy}{\mathbf{d_y}}   
\newcommand{\f}{\mathbf{f}}
\newcommand{\TT}{\mathbf{T}}
\newcommand{\kk}{\mathbf{k}}
\newcommand{\FF}{\mathbf{F}}
\newcommand{\phiphi}{\pmb{\phi}} 
\newcommand{\GammaGamma}{\pmb{\Gamma}} 
\newcommand{\XiXi}{\pmb{\Xi}} 
\newcommand{\dd}{\mathbf{d}}
\newcommand{\h}{\mathbf{h}}
\newcommand{\x}{\mathbf{x}} 
\newcommand{\widex}{\mathbf{\widetilde{x}}}
\newcommand{\widey}{\mathbf{\widetilde{y}}}
\newcommand{\wider}{\mathbf{\widetilde{r}}}
\newcommand{\y}{\mathbf{y}}
\newcommand{\M}{\mathbf{M}}
\newcommand{\I}{\mathbf{I}}
\newcommand{\II}{\mathbf{I}}
\newcommand{\PP}{\mathbf{P}}
\newcommand{\HH}{\mathbf{H}}
\newcommand{\GG}{\mathbf{G}}
\newcommand{\D}{\mathbf{D}}
\newcommand{\B}{\mathbf{B}}
\newcommand{\NN}{\mathbf{N}}
\newcommand{\LL}{\mathbf{L}}
\newcommand{\Q}{\mathbf{Q}}
\newcommand{\biggamma}{\mathbf{\Gamma}} 
\newcommand{\rr}{\mathbf{r}}
\newcommand{\capx}{\mathbf{X}}
\newcommand{\capr}{\mathbf{R}}
\newcommand{\capf}{\mathbf{F}}
\newcommand{\capp}{\mathbf{P}}
\newcommand{\capu}{\mathbf{U}}
\newcommand{\capv}{\mathbf{V}}
\newcommand{\capz}{\mathbf{Z}}
\newcommand{\capzdot}{\mathbf{\dot{Z}}}
\newcommand{\capy}{\mathbf{Y}}
\newcommand{\capxdot}{\dot{\mathbf{X}}} 
\newcommand{\cc}{\mathbf{c}}
\newcommand{\cnot}{\cc_0}
\newcommand{\pp}{\mathbf{p}} 
\newcommand{\dee}{\mathbf{d}}
\newcommand{\zero}{\mathbf{0}}
\newcommand{\batx}{\mathbf{\hat{x}}}
\newcommand{\xt}{\mathbf{x}_t}
\newcommand{\yt}{\mathbf{y}_t}
\newcommand{\rt}{\mathbf{r}_t}
\newcommand{\ex}{\mathbf{e_x}}
\newcommand{\ey}{\mathbf{e_y}}
\newcommand{\ef}{\mathbf{e_f}}
\newcommand{\eg}{\mathbf{e_g}}
\newcommand{\Jf}{\mathbf{Jf}}
\newcommand{\JJ}{\mathbf{J}}
\newcommand{\JLf}{\mathbf{J}_{\text{L}} \mathbf{f}}
\newcommand{\xset}{\mathcal{X}}
\newcommand{\zset}{\mathcal{Z}} 
\newcommand{\LS}{LS} 
\newcommand{\reg}{\mathcal{Z}^{\text{R}}}
\newcommand{\irreg}{\text{IRREGSOL}}
\newcommand{\bomega}{\pmb{\omega}}  

\newenvironment{bsmallmatrix}{\left[\begin{smallmatrix}}{\end{smallmatrix}\right]}

\title{Identifiability and Observability of Nonsmooth Systems via Taylor-like Approximations}
\author{Peter Stechlinski, Sameh Eisa, and Hesham Abdelfattah   }
\thanks{Peter Stechlinski is an associate professor in 
the department of Mathematics and Statistics, University of Maine, Maine, USA (e-mail: peter.stechlinski@maine.edu).}
\thanks{Sameh Eisa is an assistant professor in the department of Aerospace Engineering and Engineering Mechanics, University of Cincinnati, Ohio, USA (e-mail: eisash@ucmail.uc.edu).}
\thanks{Hesham Abdelfattah is a PhD student in the department of Aerospace Engineering and Engineering Mechanics, University of Cincinnati, Ohio, USA (e-mail: abdelfhm@mail.uc.edu).}

\maketitle

\begin{abstract}
New sensitivity-based methods are developed for determining identifiability and observability of nonsmooth input-output systems. More specifically, lexicographic calculus is used to construct nonsmooth sensitivity rank condition (SERC) tests, which we call lexicographic SERC (L-SERC) tests. The introduced L-SERC tests are: (i) practically implementable and amenable to large-scale problems; (ii) accurate since they directly treat the nonsmoothness while avoiding, e.g., smoothing approximations; and (iii) analogous to (and indeed recover) their smooth counterparts. To accomplish this, a first-order Taylor-like approximation theory is developed using lexicographic differentiation to directly treat nonsmooth functions. A practically implementable algorithm is proposed that determines partial structural identifiability or observability, a useful characterization in the nonsmooth setting. Lastly, the theory is illustrated through an application in climate modeling.
\end{abstract}


\section{Introduction}\label{sec.introduction}
\IEEEPARstart{I}{dentifiability} and observability are important properties of dynamical and control systems as they characterize how much information one can obtain about the parameters and states of the system, respectively, based on measurements of the system's output. 
Both properties, identifiability and observability, are longstanding problems that have received considerable attention in the literature of system identification and control theory \cite{bellman1970structural,hermann1977nonlinear,reid1977structural}.
Consider the following system:
\begin{subequations}\label{eq:1}
\begin{align}
 \dot{\x}(t) &= \f(\x(t),\uu(t),\bm{\theta}), \quad  \x(t_0) = \f_0(\btheta),\label{eq:1a}\\ 
 \y(t) &= \h(\x(t),\uu(t),\bm{\theta}),\label{eq:1b} 
 \end{align}
  \end{subequations}
where $\x(t) \in \mathbb{R}^{n_x}$ is the vector of states, $\uu(t) \in \mathbb{R}^{n_u}$ is the vector of admissible control inputs, $\bm{\theta} \in \mathbb{R}^{n_p}$ is the vector of parameters, $\y(t) \in \mathbb{R}^{n_y}$ is the vector of output variables, and $t$ is the independent variable. The right-hand side (RHS) functions $\f:D_x \times D_u \times \Theta \to \real^{n_x}$ and   $\h:D_x \times D_u \times \Theta \to \real^{n_y}$ are defined on the open and connected sets $D_x \subseteq \real^{n_x}$, $D_u \subseteq \real^{n_u}$, $\Theta \subseteq \real^{n_p}$, and $\f_0: \Theta \to D_x$ is the initial state function (for example, $\f_0(\btheta)=\x_0$ for a typical initial value problem).

In this paper, we are concerned with notions of identifiability and observability that are local in nature. First, we are concerned with local structural identifiability; the system in \eqref{eq:1} is said to be locally structurally identifiable \cite{glover1974parametrizations,dotsch1996test,stigter2015fast}
if the system parameters can be uniquely identified from the outputs in a neighborhood of some reference parameters of interest $\btheta^* \in \Theta$.
Second, we consider local observability in the sense of \cite{stigter2018efficient,van2022sensitivity}, which treats 
local observability as a special case of local identifiability; if the parameters are set to be the initial  values of the state variables $\x$, then the above-mentioned local identifiability problem becomes a local observability one since the question becomes about identifying the initial state based on the outputs. 
The reader may refer to \cite{miao2011identifiability} for more  on the problem of structural identifiability, and \cite{sontag1984concept} for more on the concept of local observability.

Various methods have been proposed in the literature for studying local structural identifiability and observability. However, many standard methods have been shown challenging to implement, especially for  large-scale/high-dimensional systems. For example, methods based on computing symbolic/analytic expressions can become cumbersome and therefore not computationally relevant  \cite{van2022sensitivity}; geometric-based methods, such as the generating series approach \cite{tunali1987new} or the observability rank condition test \cite{hermann1977nonlinear}, require successive Lie derivative computations involving the function $\h$ in \eqref{eq:1b} that are typically computationally expensive  \cite{stigter2015fast,van2022sensitivity}.
Moreover, said geometric-based methods require the system in \eqref{eq:1} to be in a control-affine form (affine in the control inputs $\uu$), which is a restrictive condition for some applications.
On the other hand, sensitivity-based methods, in particular the ``sensitivity rank condition" (SERC) \cite{dotsch1996test,stigter2015fast,stigter2018efficient,van2022sensitivity,miao2011identifiability,jacquez1985numerical}
is a computationally-relevant method for characterizing local identifiability and observability that does not suffer from the shortcomings outlined above; as has been demonstrated in \cite{van2022sensitivity}, the SERC method typically outperforms geometric-based methods computationally when applied to large-scale systems and there is no restriction that the system in \eqref{eq:1} has to be in control-affine form.

In all of the above-mentioned literature  of local structural identifiability and local observability , the RHS functions appearing in \eqref{eq:1} must be  \textit{smooth} (i.e., differentiable) to facilitate the geometric-based and sensitivity-based tools (including SERC). 
On the other hand, 
many problems 
from a wide variety of domains are concerned with physical phenomena and mathematical formulations that are inherently nonsmooth/discontinuous  \cite{Goebel2009,Cortes2012,bernardo2008piecewise}.
But  allowing for nonsmoothness in the RHS functions in \eqref{eq:1} invalidates said geometric-based and sensitivity-based methods for analyzing local identifiability and observability. Furthermore, arguments that such nondifferentiability may be ignored because it will not be reached in practice fall short even for almost-everywhere differentiable functions, for two reasons: 
(i) the real number line is represented finitely in computer simulations, and hence points of nonsmoothness may be visited with probability greater than zero  (see the diverse examples from modeling multistream heat exchangers \cite[section 4]{ref:OMS_generalizedderivatives} and investigating nonsmooth activation functions in deep learning methods  \cite{bertoin2021numerical}); and 
(ii) compositions of nonsmooth functions can lead to nondifferentiability with certainty due to the image of an inner function residing in the domain of nondifferentiability of an outer function.
Beyond ignoring the nonsmoothness, smoothing approximations prove to be popular in the literature, but this approach typically introduces  numerical error, requires user-defined inputs that are subjective, and commonly increases the complexity of the problem in  numerical treatments. Indeed, they have been shown to provide less informative analysis (see, e.g., \cite{eisa2021sensitivity}),  inaccurate predictions (see, e.g., \cite{budd2022dynamic}), and generally do not come with guarantees of convergence to the ``true'' problem (i.e., the nonsmooth version). Hence, smoothing should be avoided if computationally-relevant direct treatments of the nonsmoothness are possible. 

Based on the discussion above, we are motivated to introduce sensitivity-analysis based local identifiability and local observability methods that allow for (and directly treat) \textit{nonsmoothness} in $\f$, $\h$, and $\f_0$ in \eqref{eq:1}. We accomplish this by creating a new Taylor-like first-order approximation theory for nonsmooth functions and combining it with recent developments in nonsmooth dynamical systems theory \cite{khan2014}, based on lexicographic differentiation \cite{Nesterov2005} and lexicographic directional differentiation \cite{Khan2015_automaticdiff} (see  \cite{ref:OMS_generalizedderivatives} for an overview).  This theory is relevant to nonsmooth but continuous functions that are lexicographically smooth (L-smooth) in the sense of \cite{Nesterov2005}. Roughly, an L-smooth function is locally Lipschitz and directionally differentiable to arbitrary order. The class of L-smooth functions includes all $C^1$, convex, piecewise smooth continuous (in the sense of   \cite{bernardo2008piecewise}), and piecewise differentiable (in the sense of \cite{Scholtes}) functions, and any compositions of such functions. 
Examples of models with L-smooth functions are widespread in science (e.g., large-scale brain models \cite{coombes2018networks}, infectious disease models \cite{wang2019multiple}, climate models \cite{budd2022dynamic}, atmospheric chemistry \cite{landry2009solving} and deep learning methods \cite{liu2021modeling,bertoin2021numerical}) and engineering (e.g., power systems \cite{eisa2021sensitivity}, process engineering \cite{2017_NonsmoothDAE_CACE}, electrical circuits \cite{acary2010nonsmoothcircuits}, mechanical systems \cite{brogliato1999nonsmoothmechanics}, and bioreactors \cite{2016_Hoffner}). 
Extending identifiability/observability to this class of modeling frameworks is therefore desirable.  

 
In this paper, we introduce a novel  sensitivity-based test for local identifiability and observability of nonsmooth ODE systems in the form of \eqref{eq:1},  referred to as the lexicographic SERC (L-SERC) test. Said test is analogous to the SERC test for smooth systems, and recovers the SERC test whenever RHS functions are smooth (without needing to complete any checks for smoothness). In particular, we introduce verifiable sufficient conditions for a newly defined concept in the nonsmooth setting called ``partial'' local identifiability/observability, and we propose a practically implementable algorithm to characterize said properties in a useful way.  We also discuss the possibility of false positives and false negatives in the SERC and L-SERC tests, providing an explanation for why the smooth SERC method should only be considered a sufficient condition  and solidifying the observations in  \cite{van2022extending} on what those authors call the ``singular'' case. We apply the proposed methods to a nonsmooth climate model from \cite{ashwin2012tipping}. To accomplish the above, we create a first-order ``nonsmooth linearization'' in a Taylor-like fashion for L-smooth functions, in the spirit of classical linear approximations for smooth functions, which we hope to find widespread use elsewhere.  
The rest of the paper is structured as follows: a new first-order Taylor-like approximation theory using lexicographic derivatives is established  in Section \ref{sec:nonsmooth}. Section \ref{sec:nonsmoothSERC} presents the new sensitivity-based tools, namely, the creation of the L-SERC test. Next,  a discussion about false positives and false negatives of SERC/L-SERC tests appears in Section \ref{sec:practicalalgo}, followed by the presentation of a practical algorithm. The results are illustrated via a climate model application in Section \ref{sec:stommel} before concluding remarks are given in Section \ref{sec:conclusion}.

\section{Taylor-like Linear Approximation Theory for Nonsmooth Functions}\label{sec:nonsmooth}

Before proceeding to the new approximation theory for nonsmooth functions, we first  review generalized derivatives theory, including lexicographic differentiation. We adopt the notational convention that a well-defined vertical block matrix $\begin{bsmallmatrix} \M \\ \NN \end{bsmallmatrix}$ can be equivalently written as $(\M,\NN)$ (mimicking that a column vector $\begin{bsmallmatrix} \x \\ \y \end{bsmallmatrix}$ can be equivalently written as $(\x,\y)$).  The $n \times n$ identity matrix is denoted by $\II_n$ and the $m \times n$ zero matrix by $\zero_{m \times n}$ (with $\zero_n$ denoting the zero vector in $\real^n$). The Jacobian matrix of a $C^1$ function $\f$ at $\x$ is denoted $\JJ \f(\x)$.

 
Let $X\subseteq \real^n$ be open and $\f:X\to\real^m$ be locally Lipschitz on $X$. Clarke's generalized derivative  \cite{Clarkebook} of $\f$ at $\xnot \in X$ is the convex hull of its limiting subdifferential (also called Bouligand subdifferential), i.e.,
\begin{equation}\label{eq:Clarke}
\partial_{\rm C}\f(\xnot):=\conv \; \partial_{\rm B}\f(\xnot) ,
\end{equation}
where the limiting subdifferential is given by 
\begin{equation*} 
\partial_{\text{B}} \f(\x)
:=\left\{\FF: \exists \x_i \to \x \text{ s.t. } \x_i \in X \setminus Z_{\f}, \JJ \f(\x_{i}) \to \FF \right\},
\end{equation*}   
with $Z_{\f}$ the zero measure subset on which $\f$ is not differentiable (Rademacher's Theorem). 
Note that $\partial_{\rm C}\f(\xnot)$ is nonempty, compact and convex \cite{Clarkebook}, 
and, if $\f$ is $C^1$ at $\xnot$, then $\partial_{\rm B} \f(\xnot) = \partial_{\rm C} \f(\xnot)=\{\Jf(\xnot)\}$.
For example, 
\begin{equation*}
 f(x)=\max(0,1-e^{-x})  \Longrightarrow   \partial_{\rm C} f(x)=
  \begin{cases}
  \{0\},& x < 0, \\
  [0,1],& x = 0, \\
  \{1\},& x > 0. 
  \end{cases}
\end{equation*}
Clarke's generalized derivative is widely popular because of its usefulness in the nonsmooth numerical methods mentioned in Section \ref{sec.introduction}. 
However, it is difficult to calculate elements of Clarke's generalized derivative in general for complex functions/models because, in part, it does not satisfy sharp calculus rules. For example, consider $f(x)=\max(0,1-e^{-x})$ and $g(x)=\min(0,1-e^{-x})$. Then $h(x)=f(x)+g(x)=1-e^{-x}$, so that $\partial_{\rm C} h(0)=\{h'(0)\}=\{1\}$. But $0 \in \partial_{\rm C}f(0)$ and $0\in\partial_{\rm C}g(0)$, so that $\partial_{\rm C} h(0) \neq \partial_{\rm C}f(0) + \partial_{\rm C}g(0)$.

However, lexicographic differentiation \cite{Nesterov2005,Khan2015_automaticdiff} can be used to overcome these limitations. 
A locally Lipschitz function  $\f:X \rightarrow \mathbb{R}^m$ is lexicographically smooth (L-smooth) at $\xnot \in X$ \cite{Nesterov2005} if, for any $k \in \posint$ and directions matrix 
$$\M=[\mathbf{m}_1 \quad \mathbf{m}_2 \quad \cdots \quad \mathbf{m}_k] \in \real^{n \times k},$$ 
the homogenization sequence 
\begin{equation}\label{eq:Lsmoothfunctions}
\begin{aligned}
\f^{(0)}=\f^{(0)}_{\xnot,\bm{\mathrm{M}}} &: \real^n \to \real^m: \dd \mapsto \f'(\xnot;\dd),\\
\f^{(1)}=\f^{(1)}_{\xnot,\bm{\mathrm{M}}} &: \real^n \to \real^m: \dd \mapsto [\f^{(0)}_{\xnot,\bm{\mathrm{M}}}]'(\mathbf{m}_1;\dd),\\
\vdots\\
\f^{(k)}=\f^{(k)}_{\xnot,\bm{\mathrm{M}}} &: \real^n \to \real^m: \dd \mapsto [\f^{(k-1)}_{\xnot,\bm{\mathrm{M}}}]'(\mathbf{m}_k;\dd),
\end{aligned}
\end{equation}
exists, where $\f'(\xnot;\dd)=\lim_{\alpha \to 0^+} \frac{\f(\xnot+\alpha \dd)-\f(x)}{\alpha}$. 
That is, the class of L-smooth functions are locally Lipschitz and directionally differentiable to arbitrary order, which includes all $C^1$, $PC^1$ \cite{Scholtes}, and convex functions, and compositions thereof. 
Hence, this class of functions includes nonsmooth elementary functions such as abs-value, min, max, mid, p-norms, etc.
Roughly, the columns of the directions matrix $\M$ act as probing directions in the homogenization process (with decreasing ``importance'' from left to right), which acts to ``zoom in'' on the center and ``flatten'' the function at each stage. In line with this idea, if the columns of $\M$ span the domain space $\real^n$, it is guaranteed that the final mapping $\f^{(k)}_{\xnot,\bm{\mathrm{M}}}$ is linear \cite{Nesterov2005}.  
Nesterov defined the lexicographic derivative (L-derivative) as the derivative of this mapping:
\begin{equation}\label{eq:Lderivative}
\JJ_{\rm L}\f(\xnot;\M):=\JJ \f^{(k)}(\zero_n)\in \real^{m \times n},
\end{equation}
and  lexicographic subdifferential (L-subdifferential) as 
\begin{equation*}
\partial_{\rm L}\f(\xnot):=\{\JJ_{\rm L}\f(\xnot;\M): \M\in\real^{n \times k} \text{ has full row rank}\}.
\end{equation*}
(Note the choice of $\zero_n$ as the argument of $\JJ\f^{(k)}$ above is immaterial, because $\f^{(k)}$ is linear.)
The relevancy of the L-derivative, which is a Jacobian-like object, is summarized as follows (adapted from \cite{Nesterov2005,khan2014generalized,Khan2015_automaticdiff}): If $\f$ is L-smooth at $\xnot$, then $\partial_{\rm L}\f(\xnot) \subseteq \partial_{\rm P}\f(\xnot)$ (the plenary hull \cite{1977Sweetser} of Clarke's generalized derivative).  If $f$ is L-smooth at $\xnot$ and scalar-valued, then $\partial_{\rm L}f(\xnot) \subseteq \partial_{\rm C}f(\xnot)$. If $\f$ is $PC^1$ at $\xnot$, then $\partial_{\rm L}\f(\xnot)=\partial_{\rm B}[\f'(\xnot;\cdot)](\zero_n) \subseteq \partial_{\rm B}\f(\xnot)$.  Thanks to these relations, L-derivatives can be supplied to the nonsmooth equation-solving and optimization methods mentioned in Section \ref{sec.introduction}, and are thus computationally relevant.

Then, as a convenient tool to calculate L-derivatives, Khan and Barton \cite{Khan2015_automaticdiff}  introduced the lexicographic directional derivative (LD-derivative) of L-smooth $\f:X \rightarrow \mathbb{R}^m$ at $\xnot \in X$ in the directions $\M$ by  
\begin{equation}\label{eq:LDderivative}
\f'(\xnot;\M):=[\f^{(0)}(\mathbf{m}_1)
\;\;
\f^{(1)}(\mathbf{m}_2)
\;
\cdots
\;
\f^{(k-1)}(\mathbf{m}_k)].
\end{equation}  
If $\M$ has full row rank (and is hence right-invertible), then an LD-derivative is related to (and can be used to furnish) an L-derivative via the following linear equation system:
\begin{equation}\label{eq:LDderivative:vs:Lderivative}
\underbrace{\f'(\xnot;\M)}_{m \times k}=\underbrace{\JJ_{\rm L}\f(\xnot;\M)}_{m \times n} \underbrace{\M}_{n \times k}.
\end{equation}
If $\f$ is $C^1$, then $\JJ_{\rm L}\f(\xnot;\M)=\JJ \f(\xnot)$ and  so 
\begin{equation}\label{eq:LDderivative:vs:Lderivative:C1}
\f'(\xnot;\M)=\JJ\f(\xnot)\M.
\end{equation} 
Moreover, the LD-derivative satisfies sharp calculus rules, with closed-form expressions available for nonsmooth elemental functions such as $\max$, $\min$, absolute-value, etc. (see \cite{Khan2015_automaticdiff,ref:OMS_generalizedderivatives}). For example,   given a pair of L-smooth functions, $\f$ and $\mathbf{g}$, with appropriate domains and ranges, the LD-derivative chain rule reads 
\begin{equation}\label{eq.LDderiv.chainrule}
[\mathbf{g}\circ \f]'(\xnot;\bm{\mathrm{M}})=\mathbf{g}'(\f(\xnot);\f'(\xnot;\bm{\mathrm{M}})).
\end{equation}
The LD-derivative of the abs-value function  along $\bm{\mathrm{M}}=[m_1 \quad m_2 \quad \cdots \quad m_k] \in \mathbb{R}^{1\times k}$ is given by
\begin{equation}\label{eq.LDderivative.abs}
\mathrm{abs}'(x_0;\bm{\mathrm{M}})=\mathrm{fsign}(x_0,m_1,m_2,\dots,m_k)\bm{\mathrm{M}},
\end{equation}
  where first-sign function $\text{fsign}(\cdot)$ returns the sign of the first nonzero element, or
  zero if its input is zero.  
For familiarity, note that $\sin'(x_0;\M)=\cos(x_0)\M$, $\exp'(x_0;\M)=\exp(x_0)\M$, etc.  
The LD-derivative of the  maximum function along $\bm{\mathrm{M}}=\begin{bsmallmatrix} \M_1 \\ \M_2 \end{bsmallmatrix}\in \mathbb{R}^{2\times k}$ is given by
\begin{align}
    &\max{'}(x_0,y_0;\bm{\mathrm{M}}) 
    =\max{'}(\begin{bsmallmatrix} x_0 \\ y_0 \end{bsmallmatrix}; \begin{bsmallmatrix} \M_1 \\ \M_2 \end{bsmallmatrix})\notag\\
                &= {\rm \bf slmax}([x_0 \quad \M_1],[y_0 \quad \M_2]) \notag\\
                &:=  \begin{cases}
                          \M_1,& \text{if } \mathrm{fsign}(\begin{bsmallmatrix} x_0 \\ \M_1^{\rm T}\end{bsmallmatrix}) \leq \mathrm{fsign}(
                          \begin{bsmallmatrix} y_0\\ \M_2^{\rm T} \end{bsmallmatrix})
                          ,\\
                          \M_2,& \text{otherwise},
                      \end{cases} \label{eq.LDderivative.max}
\end{align}
i.e., the shifted-lexicographic-maximum ${\bf slmax}(\cdot)$ returns the lexicographic maximum of the two vector arguments, left-shifted by one element. 

Thanks to these sharp calculus rules, LD-derivatives can be computed in an accurate and automatable way (e.g., thanks to automatic differentiation \cite{Khan2015_automaticdiff,khan2018branch}). 
Moreover, said generalized derivatives can be supplied to dedicated nonsmooth equation-solving methods \cite{Facchinei2014,Qi1993} and optimization methods \cite{bagirov2020numerical_BGK+20}.
Hence the LD-derivative is the desired computationally-relevant tool we seek, since it can be used to furnish an L-derivative from solving a system of linear equations in \eqref{eq:LDderivative:vs:Lderivative} if $\M$ has full row rank (and thus a right inverse). Note that whenever possible $\M=\II_n$ should be chosen for simplification, but depending on the problem, the directions matrix may be passed along from an inner problem and hence may not be chosen by the user.


We are now in a position to  establish a local Taylor-like first-order approximations of L-smooth functions in the spirit of classical Taylor theory (indeed, classical linear approximations are recovered as a special case). 

\begin{theorem}\label{thm.firstorder.LDderivative}
 Let $X \subseteq \real^n$ be open and $\f:X \to \real^n$  L-smooth at $\xnot \in X$. Then, 
\begin{equation}\label{eq:firstorderapprox}
\lim_{\dd \to \zero_n} \tfrac{1}{\|\dd\|} \|\f(\xnot+\dd)-\f(\xnot)-\f'(\xnot;[\dd \quad \II_n]) \begin{bsmallmatrix} 0 \\ \dd \end{bsmallmatrix}\| 
=0.
\end{equation} 
\end{theorem}

\begin{pf}
See Appendix \ref{sec.proof1}. 
\end{pf}

We may also give a version of Theorem \ref{thm.firstorder.LDderivative} with an L-derivative instead of LD-derivative (that will be especially useful later when we return to the identifiability/observability problem), which follows immediately from Theorem \ref{thm.firstorder.LDderivative} and the fact that, for any $\dd \in \real^n$,
\begin{align}
\f'(\xnot;[\dd \quad \II_n]) \begin{bsmallmatrix} 0 \\ \dd \end{bsmallmatrix} 
&=
\f'(\xnot;[\dd \quad \II_n]) \begin{bsmallmatrix} \zero_{1 \times n} \\ \II_n \end{bsmallmatrix} \dd \notag\\
&=
\f'(\xnot;[\dd \quad \II_n]) [\dd \quad \II_n]^{-1}\dd \notag\\
&=\JJ_{\rm L}\f(\xnot;[\dd \quad \II_n])\dd, \label{eq:firstorder.LDderivative.Lderivative}
\end{align}
following \eqref{eq:LDderivative:vs:Lderivative} because $\M=[\dd \quad \II_n]$ is (right) invertible.



\begin{corollary}\label{thm.firstorder.Lderivative}
 Let $X \subseteq \real^n$ be open and $\f:X \to \real^n$  L-smooth at $\xnot \in X$. Then, 
\begin{equation}\label{eq:firstorderapprox}
\lim_{\dd \to \zero_n}  \tfrac{1}{\|\dd\|} \|\f(\xnot+\dd)-\f(\xnot)-\JJ_{\rm L}\f(\xnot;[\dd \quad \II_n]) \dd\|=0.
\end{equation} 
\end{corollary}

\begin{remark}\label{remark.firstorder.LDderivative.Lderivative}
Theorem \ref{thm.firstorder.LDderivative} and Corollary \ref{thm.firstorder.Lderivative} tell us that 
\begin{align}
\f(\x)
&\approx \f(\xnot)+\f'(\xnot;[\x-\xnot \quad \II_n]) \begin{bsmallmatrix} 0 \\ \x-\xnot \end{bsmallmatrix} \label{eq:approx:LDderivative2}\\
&=\f(\xnot)+\JJ_{\rm L}\f(\xnot;[\x-\xnot \quad \II_n]) \; (\x-\xnot)\label{eq:approx:Lderivative2},
\end{align}   
for $\x$ near $\xnot$. Moreover, in the case that $\f$ is $C^1$, observe that \eqref{eq:approx:LDderivative2}  simplifies according to \eqref{eq:LDderivative:vs:Lderivative:C1} as follows:
\begin{align*}
\f(\x)
&\approx \f(\xnot)+\f'(\xnot;[\x-\xnot \quad \II_n]) \begin{bsmallmatrix} 0 \\ \x-\xnot \end{bsmallmatrix}\\
&=\f(\xnot)+\JJ \f(\xnot) \; [\x-\xnot \quad \II_n] \; \begin{bsmallmatrix} 0 \\ \x-\xnot \end{bsmallmatrix}\\
&=\f(\xnot)+\JJ \f(\xnot) \; (\x-\xnot),
\end{align*}
i.e., the classical linear approximation is recovered.
\end{remark}

\begin{example}\label{ex:firstorder}
Consider $f(x,y)=|x^2-y^2|$ and $(x_0,y_0)=(1,1)$. Letting 
 $\M=\begin{bsmallmatrix} x-1 & 1 & 0 \\ y-1 & 0 & 1 \end{bsmallmatrix}$ 
and $z=x^2-y^2$, we note that, by \eqref{eq:LDderivative:vs:Lderivative:C1}, 
$$z'(x_0,y_0;\M)
=\JJ z(x_0,y_0) \M
=[2(x-1)-2(y-1) \; 2 \; -2].$$
And so, by the LD-derivative  rules in \eqref{eq.LDderiv.chainrule}  and \eqref{eq.LDderivative.abs},
\begin{align*}
f'(x_0,y_0;\M)
&=
{\rm abs}^{'}(z(x_0,y_0);z'(x_0,y_0;\M))\\
&=\fsign(0,2(x-1)-2(y-1),2,-2) \\
&\quad \times [2(x-1)-2(y-1) \quad 2 \quad -2].
\end{align*}
Hence, according to \eqref{eq:approx:LDderivative2},
\begin{align*}
f(\x) 
&\approx f(\x_0)+f'(\x_0;[\x-\x_0 \quad \II_n]) \begin{bsmallmatrix} 0  \\ \x-\x_0 \end{bsmallmatrix}\\
&=\fsign(2(x-1)-2(y-1),2,-2) \\
&\quad \times (2(x-1)-2(y-1))\\
&=\begin{cases} 
2(x-1)-2(y-1) & \text{if } x \geq y, \\ 
-2(x-1)+2(y-1) & \text{if } x<y, \end{cases}
\end{align*}
which is indeed a first-order approximation near $\xnot$. Note that  $\dd \mapsto \JJ_{\rm L} \f(\xnot;[\dd \quad \II_n])$ (and $\dd \mapsto \f'(\xnot;[\dd \quad \II_n])$) need not be continuous in general, as demonstrated in this example where 
 \begin{align*}
 J_{\rm L} f(\xnot;[\dd \quad \II])
 &=\begin{cases} 1, &\text{if } d_1 \geq d_2,\\ -1, &\text{if } d_1<d_2.\end{cases}
 \end{align*} 
 \end{example}


The Taylor-like approximation theory established above builds on directional derivative first-order approximations \cite{Scholtes}. Recently,  Griewank \cite{griewank2013stable} (see also \cite{Griewank2020beyond}) established an abs-linearization theory (a first-order nonsmooth linearization)   for functions in abs-normal form, i.e., those composed of finitely many smooth elemental functions and the abs-value function. Streubel et al.\ \cite{streubel2021piecewise} and Griewank et al.\ \cite{griewank2021abs}   expanded on this technique to provide a generalized Taylor expansion of arbitrary order for abs-normal functions, again mirroring the classical Taylor expansion of smooth functions.  This theory is computationally relevant (it is motivated by, and built with, automatic differentiation in mind), allowing for efficient and accurate numerical implementations (see the discussion at the bottom of \cite[Page 341]{Griewank2020beyond}). 
Moreover, this approach has the added benefit of being ``semi-local'' in the sense that a generalized expansion at a smooth domain point can encode nonsmooth  behavior away from said point (see  the half-pipe example on \cite[Page 341]{Griewank2020beyond}). However, this powerful Taylor-like toolkit for nonsmooth functions is not suitable for the current problem of interest (sensitivity-based analysis of identifiability and observability) for two reasons: (i) solutions of input-output nonsmooth ODEs with abs-normal form RHS functions need not admit solutions that are abs-normal with respect to parameters (see the counterexample in \cite[Page 1056]{ref:OMS_generalizedderivatives}); and (ii) the form of the abs-linearization in those works is not amenable to the ordinary least squares derivations involved in SERC and L-SERC (just like  directional derivative first-order approximations are not --- see the comments in Remark \ref{re.LSERC}). 
The flexibility of the LD-derivative Taylor-like approximation theory in this work overcomes these issues, while remaining  computationally relevant.

\section{Nonsmooth Sensitivity Rank Condition for Local Identifiability and Observability}\label{sec:nonsmoothSERC}

In this part, we develop theory for local identifiability and observability in the nonsmooth setting via an LD-derivative SERC-like test. First, we provide the definition of local structural identifiability  \cite{stigter2015fast,miao2011identifiability} for the nonlinear system in  \eqref{eq:1}, under the following assumptions:

\begin{assumption}\label{assumption.Lsmooth}
Suppose that the RHS functions $\f$, $\f_0$, and $\h$ in \eqref{eq:1} are L-smooth on their respective domains. Suppose that the set of admissible controls is given by $\mathcal{U}:=L^1([t_0,t_f],D_u)$, i.e., Lebesgue integrable. 
\end{assumption}

\begin{assumption}\label{ass.solution}
Suppose that the system in \eqref{eq:1} admits an (absolutely continuous)  solution $\x^*(t):=\x(t;\uu^*,\btheta^*)$ on $[t_0,t_f]$, for the given reference parameters $\bm{\theta^}* \in \bm{\Theta}\subseteq \mathbb{R}^{n_p}$ and reference control input $\uu^* \in \mathcal{U}$.
\end{assumption}
 
Note that L-smoothness (and thus local Lipschitz continuity) of the RHS functions in Assumption \ref{assumption.Lsmooth} implies the existence of a solution to \eqref{eq:1} but only locally, hence the need for Assumption \ref{ass.solution} which is typical in the smooth setting too. Note also that uniqueness of $\x^*$ follows by local Lipschitz continuity of  $(\x,\uu,\btheta) \mapsto \f(\x,\uu,\btheta)$.
 
\begin{definition}\label{defn:identifiability}
The system in \eqref{eq:1} is  \textit{locally structurally identifiable} at $\bm{\theta^}*$ if, for any  $\bm{\theta}_1$ and $\bm{\theta}_2$ within a  neighborhood of $\bm{\theta^}*$, we have that  $\y(t;\uu^*,\bm{\theta_1})=\y(t;\uu^*,\bm{\theta_2)}$ for all $t\in[t_0,t_f]$ if and only if $\bm{\theta_1}=\bm{\theta_2}$.
\end{definition}
  
It is worth mentioning that the term ``structural" has been used to refer to different concepts by different authors \cite{distefano1980parameter}. However, in this paper, similar to \cite{stigter2015fast}, we consider the term ``structural" to mean that the local identifiability property is determined by the model's structure, without taking into consideration the effect of noise and uncertainty in measurements, as opposed to practical identifiability \cite{wieland2021structural}, which takes into account uncertainties and noise involved in the measurements (e.g., due to sensors).  
Going forward, we will drop the use of the term ``local,'' though it is implicitly understood that these properties are local in nature for nonlinear systems.  

Before proceeding with the nonsmooth setting,  we provide context via the smooth setting: the sensitivity rank condition (SERC) \cite{dotsch1996test,stigter2015fast,stigter2018efficient,van2022sensitivity,miao2011identifiability,jacquez1985numerical}, mentioned earlier in the introduction, is a sensitivity-based analysis that can be used to determine local identifiability or observability for systems in the form \eqref{eq:1} if Assumption  \ref{assumption.Lsmooth} is replaced with the stronger requirement that the RHS functions are $C^1$. Under this assumption, the so-called forward sensitivity functions 
\begin{equation}\label{eq.sensitivityfunctions}
\caps_{\x}(t):=\frac{\partial \x}{\partial \btheta}(t;\uu^*,\btheta^*), \quad \caps_{\y}(t):=\frac{\partial \y}{\partial \btheta}(t;\uu^*,\btheta^*)
\end{equation}
of $\x$ and $\y$, respectively, 
corresponding to the reference parameters $\bm{\theta^*}$ and reference control $\uu^*$, are the unique solutions on $[t_0,t_f]$ of the forward sensitivity system
\begin{subequations}\label{eq:sens.smooth}
\begin{align}
\dot{\caps}_{\x}(t)
&= \frac{\partial \f}{\partial \x} \caps_{\x}(t)+ 
\frac{\partial \f}{\partial \bm{\theta}}, \quad \caps_{\x}(t_0)
=\JJ \f_0(\btheta^*)\label{eq:sens.smooth1}\\
\caps_{\y}(t) 
&= \frac{\partial \h}{\partial \x} \caps_{\x}(t)+ 
\frac{\partial \h}{\partial \bm{\theta}}.\label{eq:sens.smooth2}
\end{align}
\end{subequations} 
with omitted partial derivative arguments  $(\x^*(t),\uu^*(t),\bm{\theta^*})$.

The spirit of the sensitivity-based identifiability test arises from that fact that (ignoring arguments temporarily for simplicity) we can write the first-order approximation $\Delta \y=\caps_{\y} \Delta \btheta$, where $\Delta \bm{\theta} = \bm{\theta} - \bm{\theta^*}$ and $\Delta \y = \y(\bm{\theta})-\y(\bm{\theta^*})$. This motivates a rank test of $\caps_{\y}$, since invertibility yields $\Delta \btheta=\left(\caps_{\y} \right)^{-1}\Delta \y$, i.e.,  changes in the output can be mapped to unique  changes in the parameters. However, since there are typically fewer outputs than (unknown) model parameters (i.e., $n_y < n_p$), left-invertibility of $\caps_{\y}$ may not be possible, motivating the so-called ``SERC matrix,'' which is constructed by sampling the output sensitivity functions $\caps_{\y}$ at times $t_0, t_1, t_2, \ldots, t_N \in [t_0,t_f]$, with $N \in \posint$ satisfying $(N+1)n_y \geq n_p$, and stacking the results into a block matrix as follows:
\begin{align}\label{SERC_ID}
{\pmb \Upsilon}:=
(
\caps_{\y}(t_0),
\caps_{\y}(t_1),
\ldots
\caps_{\y}(t_N))  \in \mathbb{R}^{(N+1)n_y \times n_p}.
\end{align}
Then, using the first-order approximation $\Delta \y=\caps_{\y} \Delta \btheta$, at any sample $k=0,1,2,...,N$, we get
$$\y(t_k;\uu^*(t_k),\btheta) 
\approx 
\y(t_k;\uu^*(t_k),\btheta^*)+\caps_{\y}(t_k) \; (\btheta-\btheta^*),$$ 
and an ordinary least squares derivation (see, e.g., \cite[section 5]{miao2011identifiability}) can be used to show that the residual sum of squares between exact measurements and first-order approximation attains a (local) minimum at $\btheta=\btheta^*$ if ${\pmb \Upsilon}^{\rm T}{\pmb \Upsilon}$ is nonsingular. Since $\text{rank}({\pmb \Upsilon})=\text{rank}({\pmb \Upsilon}^{\rm T}{\pmb \Upsilon})$, 
the SERC identifiability test can be stated as follows \cite{stigter2015fast,stigter2018efficient,van2022sensitivity,van2022extending}, based on the notion of ``sensitivity identifiability'' (which we call SERC identifiability here): 

\begin{definition}\label{defn:identifiability.SERC}
The system in \eqref{eq:1} is said to be \textit{SERC identifiable} at $\bm{\theta^}* \in \bm{\Theta}$ if $\text{rank}({\pmb \Upsilon})=n_p$.  
\end{definition}

The relevance of this test follows from the next result, which we emphasize again requires $C^1$ smoothness of the RHS functions in \eqref{eq:1}.


\begin{theorem}\label{thm:smoothSERC.identifiability}
If the system in \eqref{eq:1} is SERC identifiable at $\btheta^*$, then it is structurally identifiable at $\btheta^*$. 
\end{theorem}

With the above in mind, we return to the nonsmooth setting and build a SERC-like test for  identifiability. To accomplish this, we use the nonsmooth sensitivity theory developed by Khan and Barton in \cite{khan2014} for ODE systems with L-smooth RHS functions, adjusted for an input-output system. Namely, under Assumptions \ref{assumption.Lsmooth}-\ref{ass.solution}, it follows that $\x^*_t:=\x^*(t;\uu^*,\cdot)$ and $\y^*_t:=\y(t;\uu^*,\cdot)$ are L-smooth functions; given any direction $\dd \in \real^{n_p}$, the LD-derivative mappings
\begin{align}
\capx^*(t)&=[\x^*_t]'(\btheta^*;[\dd \quad \II_{n_p}]) \in \real^{n_x \times (1+n_p)},\label{eq:Xsens}\\ 
\capy^*(t)&=[\y^*_t]'(\btheta^*;[\dd \quad \II_{n_p}]) \in \real^{n_y \times (1+n_p)},\label{eq:Ysens}       
\end{align}
are absolutely continuous functions that uniquely satisfy the nonsmooth forward sensitivity system   on $[t_0,t_f]$:
\begin{align}
&\dot{\capx}(t)
=\f'(\x^*(t),\uu^*(t),\btheta^*;(\capx(t),\zero_{n_u \times (1+n_p)},[\dd \;\; \II_{n_p}])),\notag\\
&\capy(t)
=\h'(\x^*(t),\uu^*(t),\btheta^*;(\capx(t),\zero_{n_u \times (1+n_p)},[\dd \;\;  \II_{n_p}])),\notag\\
&\capx(t_0)
=\f_0'(\btheta^*;[\dd \;\; \II_{n_p}]).\label{eq:sens.nonsmooth}
\end{align}  
Motivated by the SERC matrix in \eqref{SERC_ID},  we define the lexicographic SERC (L-SERC) matrix corresponding to the ``primary probing direction'' $\dd \in \real^{n_p}$  as follows:
\begin{equation}\label{LSERC_ID}
{\pmb \Upsilon}_{\dd}:=
\begin{bmatrix} 
\caps_{\y}^{\rm L}(t_0) \\
\caps_{\y}^{\rm L}(t_1) \\
\vdots \\
\caps_{\y}^{\rm L}(t_N) 
\end{bmatrix}
= 
\begin{bmatrix} 
\text{lshift}(\capy^*(t_0))   \\
\text{lshift}(\capy^*(t_1)) \\
\vdots \\
\text{lshift}(\capy^*(t_N))
\end{bmatrix} \in \real^{(N+1)n_y \times n_p},
\end{equation}
where $\text{lshift}(\cdot)$ denotes the left-shift operator and
where the output L-sensitivity functions are given by
 \begin{equation}\label{eq.Lsensitivities}
\caps_{\y}^{\rm L}(t):=\capy^*(t) [\dd \quad \II_{n_p}]^{-1} \in \partial_{\rm L}\y^*_t(\btheta^*).
 \end{equation}
Note that $\dd$ is called the  primary probing direction because the leftmost column in the directions matrix determines the ``piece'' of the nonsmooth function to probe into, and is thus the most important in that sense.
Before proceeding, we note that structural identifiability (Definition \ref{defn:identifiability}) still applies in this setting. However, we propose the following specialized definition  of identifiability, motivated by the complications arising from the presence of nonsmoothness now being considered.
  
\begin{definition}\label{defn:identifiability.partial} 
The system in \eqref{eq:1} is said to be \textit{partially structurally  identifiable} at $\bm{\theta^}* \in \bm{\Theta}$ if there exist a neighborhood $N \subseteq \bm{\Theta}$ of $\btheta^*$ and a connected set $V \subseteq \bm{\Theta}$ containing $\btheta^*$ such that  for any $\bm{\theta}_1, \bm{\theta}_2 \in N \cap V$, we have that $\y(t;\uu^*,\bm{\theta_1})=\y(t;\uu^*,\bm{\theta_2)}$  
for all $t\in[t_0,t_f]$ if and only if $\bm{\theta_1}=\bm{\theta_2}$. 
\end{definition}

With this in mind, we propose L-SERC identifiability.

 
\begin{definition}\label{defn:identifiability.partial.indirection}
The system in \eqref{eq:1} is said to be \textit{L-SERC identifiable} at $\bm{\theta^}* \in \bm{\Theta}$ in the direction $\dd \in \real^{n_p}$ if $\text{rank}({\pmb \Upsilon}_{\dd})=n_p$. 
\end{definition}

Hence, if \eqref{eq:1} is L-SERC identifiable at $\btheta^*$ in the direction $\dd$, then it is  partially structurally identifiable with $V=\{\alpha \dd:\alpha \geq 0\}$ and some neighborhood $N$ of $\btheta^*$, a subset of the   ``true'' region in parameter space in which identifiability holds.  We are now in position to present  the proposed L-SERC method.

\begin{theorem}\label{thm:nonsmoothSERC.identifiability}
Let Assumption \ref{assumption.Lsmooth} hold. If \eqref{eq:1} is L-SERC identifiable at $\btheta^*$ in the direction $\dd \in \real^{n_p}$, then it is partially structurally identifiable at $\btheta^*$. 
\end{theorem} 

\begin{pf} 
First we argue  \eqref{eq:sens.nonsmooth} is the nonsmooth sensitivity system associated with \eqref{eq:1}, and is solved by $(\capx,\capy)$ in \eqref{eq:Xsens}-\eqref{eq:Ysens}:  Note that \eqref{eq:1a}  can be rewritten as
\begin{equation} 
\dot{\x}(t) =\bar{\f}(t,\x(t),\bm{\theta}), \quad \x(t_0) = \f_0(\btheta),\label{eq:bar1}
\end{equation}
with $\bar{\f}:(t,\x,\btheta) \mapsto \f(\x,\uu^*(t),\btheta)$, and is therefore solved by $\x^*$ from Assumption \ref{ass.solution}. Moreover, thanks to Assumption \ref{assumption.Lsmooth} (and recalling that $\uu^*(t) \in D_x$ for all $t\in[t_0,t_f]$), there exists an open set $W$, satisfying 
$$X:=\{(\x^*(t),\btheta^*):t \in[t_0,t_f]\} \subset W \subseteq D_x \times \Theta,$$ 
on which: (i)  $\bar{\f}(t,\cdot,\cdot)$ is L-smooth for each $t \in [t_0,t_f]$ (from L-smoothness of $\f$); (ii) $\bar{\f}(\cdot,\x,\btheta)$ is measurable for each $(\x,\btheta) \in W$ (from Lebesgue integrability, and hence measurability, of $\uu^*$); and (iii) $\f$ is bounded on $[t_0,t_f] \times N$ (by compactness of $X$ with $W$ shrunk as needed).

Hence,  \cite[Theorem 4.2]{khan2014} (see also its restatement on \cite[Page 1055]{ref:OMS_generalizedderivatives}, combined  with the remarks regarding Equation (34) in \cite{khan2014}) may be applied to \eqref{eq:bar1}, which immediately yields that  $\x_t^*:=\x(t;\uu^*,\cdot)$ is L-smooth on a neighborhood of $\btheta^*$. Moreover, observing that  
\begin{equation}\label{eq.zerothing}
[\bar{\f}(t,\cdot,\cdot)]'(\x,\btheta;(\capx,\M))=[\bar{\f}]
'(t,\x,\btheta;(\zero,\capx,\M))
\end{equation}
for any directions matrix $\M$, it also follows from \cite[Theorem 4.2]{khan2014} that  the LD-derivative mapping $t \mapsto [\x^*_t]'(\btheta^*;\M)$ uniquely solves the first equation in \eqref{eq:sens.nonsmooth} on $[t_0,t_f]$, with initial condition given by the third equation. Lastly,  since $\h$ is L-smooth on its (open and connected) domain and $\y^*$ satisfies \eqref{eq:1b}, it follows that $\y^*_t$ is L-smooth on a neighborhood of $\btheta^*$ as a composition of L-smooth functions and, from a simple application of the LD-derivative chain rule as in \eqref{eq.zerothing}, the second equation in  \eqref{eq:sens.nonsmooth} follows.

With the nonsmooth sensitivity theory established, and Theorem \ref{thm.firstorder.LDderivative} in place, we may proceed similarly to the derivation outlined in \cite[Section 5.1]{miao2011identifiability}: given the unique solution $\x^*$, we denote the corresponding output at sampling time  $t_k \in \{t_0,t_1,\ldots,t_N\}$ by $$\y_k(\btheta):=\y(t_k;\uu^*;\btheta) \in \real^{n_y}.$$ 
Letting $\Delta \btheta=\btheta-\btheta^*$, it follows from Corollary \ref{thm.firstorder.Lderivative} (and see Remark \ref{remark.firstorder.LDderivative.Lderivative}) that
\begin{equation}\label{eq:localapprox}
\y_k(\btheta) \approx \y_k(\btheta^*)+\JJ_{\rm L}\y_k(\btheta^*;[\Delta \btheta \quad \II_{n_p}]) \Delta \btheta. 
\end{equation}
Then, letting $\rr_k$ denote the errorless measurement at $t=t_k$, the residual sum of squares (RSS) between exact and actual measurements is  
$$
\textrm{RSS}(\Delta \btheta)
 :=\sum_{k=0}^{N} ||\rr_k-\y_k(\btheta)||^2,$$
from which the approximation yields
\begin{equation*} 
\textrm{RSS}(\Delta \btheta)
\approx \sum_{k=0}^{N} ||\rr_k-\y_k(\btheta^*)-\JJ_{\rm L}\y_k(\btheta^*;[\Delta \btheta \quad \II_{n_p}]) \Delta \btheta||^2.
\end{equation*}
Then, since $\rr_k=\y_k(\btheta^*)$ by assumption, 
\begin{align}
\textrm{RSS}(\Delta \btheta)
& \approx \sum_{k=0}^{N} ||\JJ_{\rm L}\y_k(\btheta^*;[\Delta \btheta \quad \II_{n_p}]) \Delta \btheta||^2\label{eq:RSS.approx}\\
&=({\pmb \Upsilon}_{\dd} \Delta \btheta)^{\rm T}({\pmb \Upsilon}_{\dd} \Delta \btheta)=\Delta \btheta^{\rm T} {\pmb \Upsilon}_{\dd}^{\rm T} {\pmb \Upsilon}_{\dd} \Delta \btheta.\notag
\end{align} 
It follows that $\textrm{RSS}(\btheta)$ is minimized when 
\begin{equation}\label{eq:RSSequation}
{\pmb \Upsilon}_{\dd}^{\rm T} {\pmb \Upsilon}_{\dd} \Delta \btheta=\zero_{n_p}
\end{equation}
holds. Since $\text{rank}({\pmb \Upsilon}_{\dd}^{\rm T}{\pmb \Upsilon}_{\dd})=\text{rank}({\pmb \Upsilon}_{\dd})=n_p$ and $(N+1)n_y \geq n_p$, 
then \eqref{eq:RSSequation} is uniquely solved by $\Delta \btheta=\zero$, and therefore (local) identifiability holds in the direction $\Delta \btheta$. 
Lastly, observe that \eqref{LSERC_ID} follows from the fact that 
$$\caps_{\y}^{\rm L}(t_k)=\capy(t_k) [\dd \quad \II_{n_p}]^{-1}=\capy(t_k) \begin{bsmallmatrix} \zero_{1 \times n_p} \\ \II_{n_p} \end{bsmallmatrix}. \qed$$  
\end{pf}

\begin{remark}\label{re.LSERC}
Some remarks are in order:
\begin{itemize}
\item The nonsmooth sensitivity theory in \cite{khan2014} is given for non-autonomous ODE systems. Consequently, the L-SERC test in Theorem \ref{thm:nonsmoothSERC.identifiability} can be extended to the non-autonomous case of \eqref{eq:1} in a straightforward way.
\item The zero matrices appear  in the directions matrix in \eqref{eq:sens.nonsmooth} because $\uu^*$ does not depend on $\btheta$, in the same way that in the smooth case, the terms $\frac{\partial \f}{\partial \uu}$ and $\frac{\partial \h}{\partial \uu}$ do  not appear in the forward sensitivity  system in \eqref{eq:sens.smooth}. 
\item The function $\caps^{\rm L}_{\y}$  in \eqref{eq.Lsensitivities}, and the similarly defined $\caps_{\x}^{\rm L}(t):=\capx^*(t)[\dd \quad \II_{n_p}]^{-1}$, are the ``L-sensitivity functions'' of \eqref{eq:1} with respect to $\btheta^*$ and $\uu^*$ because they are nonsmooth analogs of classical (forward) sensitivity functions. Indeed, if the RHS functions in \eqref{eq:1} are $C^1$, then the classical sensitivity functions are recovered, i.e., $\caps_{\x}^{\rm L}=\frac{\partial \x}{\partial \btheta}=\caps_{\x}$ and $\caps_{\y}^{\rm L}=\frac{\partial \y}{\partial \btheta}=\caps_{\y}$.
\item Despite directional derivatives possessing a first-order approximation \cite{Scholtes}, the RSS derivation above cannot be completed using that theory. Indeed, from \cite[Theorem 3.1.2]{Scholtes},   \eqref{eq:localapprox} may be replaced with $\y_k(\btheta) \approx \y_k(\btheta^*)+\y_k'(\btheta^*;\Delta \btheta)$, 
but then \eqref{eq:RSS.approx} becomes $\textrm{RSS}(\Delta \btheta) \approx 
\sum_{k=0}^{N} \|\y_k'(\btheta^*;\Delta \btheta)\|^2$,
leading to an impasse. 
\end{itemize}
\end{remark}

 
We demonstrate how the L-SERC test can be implemented for a nonsmooth system with a simple example, which is an adaptation of \cite[Example 3.3]{ackley2021determining}. 

\begin{example}\label{ex:riot} 
Consider the following system 
\begin{align}\label{riot1}
\begin{split}
\dot{x}(t)&=\mathrm{max}(0,1-e^{-(x(t)-\theta_1)}), \quad x(0)=\theta_2,\\
y(t)&=x(t), 
\end{split}
\end{align}
where $\bm{\theta}=(\theta_1,\theta_2) \in \Theta=\mathbb{R}^2$ are the parameters. Supposing that $\btheta^*=(1,1)$ are the reference parameters, and $[t_0,t_f]=[0,t_f]$ is the time horizon of interest, we proceed as follows to construct the L-SERC matrix:  
let $f(x,\btheta)=\max(0,1-e^{-(x-\theta_1)})$, $h(x,\btheta)=x$,  and $f_0(\btheta)=\theta_2$.  Then, letting $q=q(x,\btheta)=0$ and $r=r(x,\btheta)=1-e^{-(x-\theta_1)}$, we have  $f(x,\btheta)=\max(q(x,\btheta),r(x,\btheta))$, and so, by the LD-derivative chain rule in \eqref{eq.LDderiv.chainrule} and max rule in \eqref{eq.LDderivative.max},
\begin{align*}
&f'(x,\btheta;(\capx,\M))\\
&={\max}'(q,r;(\JJ q(x,\btheta) \begin{bsmallmatrix} \capx \\ \M \end{bsmallmatrix},\JJ r(x,\btheta) \begin{bsmallmatrix} \capx \\ \M \end{bsmallmatrix}))\\
&={\rm \bf slmax}([q \quad \JJ q (x,\btheta) \begin{bsmallmatrix} \capx \\ \M \end{bsmallmatrix}],[r \quad \JJ r(x,\btheta) \begin{bsmallmatrix} \capx \\ \M \end{bsmallmatrix}])\\
&={\rm \bf slmax}([0 \quad [0 \;\; 0 \;\; 0] \begin{bsmallmatrix} \capx \\ \M \end{bsmallmatrix}],[1-z \quad [z  \;\; -z  \;\; 0] \begin{bsmallmatrix} \capx \\ \M \end{bsmallmatrix}])\\
&={\rm \bf slmax}([0 \quad \zero_{1 \times 3}],[1-z \quad z \; \capx-z \; \mathrm{row}_1(\M)]),
\end{align*}
where $z=z(x,\btheta)=e^{-(x-\theta_1)}$.
Next, by \eqref{eq:LDderivative:vs:Lderivative:C1}, 
\begin{align*}
&h'(x,\btheta;(\capx,\M))
=\JJ h(x,\btheta) \begin{bsmallmatrix} \capx \\ \M \end{bsmallmatrix}
=[1 \quad 0 \quad 0]\begin{bsmallmatrix} \capx \\ \M \end{bsmallmatrix}=\capx,\\
&f_0'(\btheta;\M)=\JJ f_0(\btheta) \M=[0 \quad 1]\M=\mathrm{row}_2(\M).
\end{align*} 
Hence, the nonsmooth sensitivity system in \eqref{eq:sens.nonsmooth} for directions matrix 
$\M=[\dd \quad \II_2]=\begin{bsmallmatrix} d_1 & 1 & 0 \\ d_2 & 0 & 1 \end{bsmallmatrix}$ 
is given by 
\begin{align}  
\dot{\capx}(t)&={\rm \bf slmax}\left(
[0 \quad \zero_{1 \times 3}], [1-z^*(t) \quad A^*(t)]
\right),\notag\\
\capy(t)&=\capx(t),\label{riot3}\\
\capx(0)&=\mathrm{row}_2(\M),\notag
\end{align}
where $z^*(t):=z(x^*(t),\btheta^*)$, 
$x^*(t)$ is the reference solution of \eqref{riot1} when $\btheta=\btheta^*$, and $A^*(t):=z^*(t) (\capx(t)-\mathrm{row}_1(\M))$. Then \eqref{riot1} and \eqref{riot3} can be numerically solved simultaneously on $[0,t_f]$ to furnish $\x^*$, $\y^*$, $\capx^*$, $\capy^*$, and therefore ${\pmb \Upsilon}_{\dd}$ in \eqref{LSERC_ID} by sampling $\capy^*(t_k)$ and  left-shifting.

In this case, we may proceed analytically  since we can calculate $x^*(t)=x(t;\btheta^*)=1$. Then, $\dd=\coord_1=(1,0)$ 
yields $\capx^*(t)=[0 \quad 0 \quad 1]=\capy^*(t)$, from which it follows that $\caps_{\y}^{\rm L}(t)=[0 \quad 1]$.
Since $n_y=1$ and $n_p=2$ here, $\pmb{\Upsilon}_{\coord_1}$ in \eqref{LSERC_ID} can be constructed by sampling at two times (the choice of sampling times is unimportant here because $\caps_{\y}$ is constant), yielding
$${\pmb \Upsilon}_{\coord_1}=\begin{bmatrix} 0 & 1 \\ 0 & 1 \end{bmatrix} \quad \Longrightarrow \quad \mathrm{rank}({\pmb \Upsilon}_{\dd})=1<n_p.$$
Hence, \eqref{riot1} is not L-SERC identifiable at $\btheta^*$ in the direction $\coord_1$.
On the other hand, with $\dd=-\coord_1$, we get $\capx^*(t)=[e^{t}-1 \quad 1-e^{t} \quad e^t]=\capy^*(t)$, so that
 $\caps_{\y}^{\rm L}(t)=[1-e^t \quad e^{-t}]$.  
Sampling at $t_0=0$, $t_1>t_0$, we get 
$${\pmb \Upsilon}_{-\coord_1}=\begin{bmatrix} 0 & 1 \\ 1-e^{t_1} & e^{t_1} \end{bmatrix} \quad \Longrightarrow \quad \mathrm{rank}({\pmb \Upsilon}_{\dd})=2=n_p.$$
That is, \eqref{riot1} is L-SERC identifiable at $\btheta^*$ in the direction $\coord_1$, and is therefore partially identifiable.

The L-SERC test results above can be verified because the ODE system in \eqref{riot1} admits a closed-form solution:
$$
x(t;\bm{\theta})=\begin{cases}\ln(e^t+e^{\theta_1-\theta_2}-e^{t+\theta_1-\theta_2}) + \theta_2, & \mathrm{if } \theta_1 \leq \theta_2, \\
\theta_2, & \mathrm{if } \theta_1 > \theta_2. \end{cases}
$$
Thus,  we can see from the closed-form solution that \eqref{riot1} is partially identifiable at $\btheta^*$  with $V=\{\btheta:\theta_1 \leq \theta_2\}$ and unidentifiable otherwise as $\theta_1$ cannot be determined from the output in that case, verifying the L-SERC test result.  
\end{example}
 
 Local observability is a special case of local identifiability 
\cite{stigter2018efficient,van2022sensitivity} since it corresponds to the case when the initial conditions act as the parameters, i.e., $\f_0(\btheta)=\btheta=\x_0$, so that $n_p=n_x$. It is therefore straightforward to establish an L-SERC test for observability:  
the system in \eqref{eq:1}  is partially observable at $\btheta^*$ if $(N+1)n_y \geq n_x$ and $\mathrm{rank}({\pmb \Upsilon_{\dd}})=n_x$, with $\capy^*$  in  \eqref{LSERC_ID} obtained instead via  
\begin{equation}\label{eq:sens.obs}
\begin{aligned}
\dot{\capx}(t)
&=\f'(\x^*(t),\uu^*(t);(\capx(t),\zero_{1 \times (1+n_u)})),\\
\capy(t)
&=\h'(\x^*(t),\uu^*(t);(\capx(t),\zero_{1 \times (1+n_u)})),\\
\capx(t_0)
&=[\dd \quad \II_{n_x}],
\end{aligned}
\end{equation} 
replacing the sensitivity system in \eqref{eq:sens.nonsmooth} for identifiability.

\section{Building a Practical L-SERC Algorithm}
\label{sec:practicalalgo}
We proceed now towards a practically implementable algorithm that characterizes useful identifiability information. We build this algorithm motivated by both ``false positive'' L-SERC tests and ``false negative'' L-SERC tests. First we note that, unsurprisingly, partial identifiability need not imply (full) structural identifiability, as demonstrated in Example \ref{ex:riot}.  Motivated by this, we should conduct several L-SERC tests with different primary probing directions (i.e., a ``primary probing stage''); we propose a ``natural'' test of identifiability (natural observability can be similarly defined), motivated by the natural (or standard/canonical) basis vectors $\{\coord_i\}$ in $\real^{n_p}$.

\begin{definition}\label{defn:identifiability.natural}
The system in \eqref{eq:1} is said to be \textit{naturally identifiable} at $\bm{\theta^}* \in \bm{\Theta}$ if $\mathrm{rank}({\pmb \Upsilon}_{\dd_i})=n_p$ for each $\dd_i=\pm \coord_i$, $i=1,\ldots,n_p$. 
\end{definition} 

Thus, \eqref{eq:1} is naturally identifiable at $\btheta^*$ if perturbations of parameters in any standard basis vector direction (positive or negative) can be uniquely identified via output measurements. 
But natural identifiability  need not imply structural identifiability either (though we suspect it to typically act as a certificate of identifiability  in practical applications). This can be illustrated with a simple single-input, single-output relation $y=|\theta|$; in this case, $\Upsilon_d=1$ for $d=e_1=1$ and $\Upsilon_d=-1$ for $d=-e_1=-1$. Hence, the L-SERC matrices have full rank in both cases, and natural identifiability holds at $\theta^*=0$ according to the definition. However, structural identifiability does not hold at $\theta^*=0$ since $\pm \theta$ yield the same output. Hence, one or several ``positive L-SERC tests'' (i.e., calculating full rankness) need not imply structural identifiability;   a ``false positive'' is possible in this sense.   However, the almost-everywhere differentiability of Lipschitzian functions does provide some relief here: if a probing direction points along a nonsmoothness, a perturbation may reveal a clue of this by returning a different rank.  Motivated by this, we consider perturbing the primary probing directions to produce a set of ``twin'' directions (i.e., a ``twin probing stage'').

Next, we consider the opposite case, i.e., ``false negative'' L-SERC tests: As noted by the authors in \cite{van2022extending}, the smooth SERC test may result in a ``false negative,'' wherein there is rank deficiency but actually the problem is identifiable. 
The same can be true in the nonsmooth setting:  consider 
 the single-input, single-output relation $y=\max(\theta^3,\theta^5)$, which is an invertible L-smooth map at $\theta^*=0$ but ``fails'' the L-SERC test for full rankness since $\Upsilon_{d}=0$ for any $d \neq 0$. 
Hence, we should proceed with caution when claiming unidentifiability of a complicated problem/model using a sensitivity-based analysis.
The authors in \cite{van2022extending} called such a $\btheta^*$ in smooth SERC an ``isolated output sensitivity singularity'' (away from $\theta^*$, this problem is avoided), and went on to propose an algorithm that tests for this degenerate case using a singular value decomposition of the SERC matrix (see also \cite{stigter2015fast,stigter2018efficient,van2022sensitivity,van2022extending}).  

With the above discussion in mind concerning false negatives, we are motivated to consider a ``singularity probing stage'' in the nonsmooth setting by analyzing the singular value decomposition (SVD) of ${\pmb \Upsilon}_{\dd}$ to help  determine if rank deficiency only holds at a single point and any perturbation yields full rankness (i.e., a false negative).  Letting $n_s:=(N+1)n_y$, and recalling that $n_s \geq n_p$, we suppose that  
\begin{equation}\label{eq:LSERC.SVD}
{\pmb \Upsilon}_{\dd}=\capu \Sigma \capv^{\rm T}=\sum_{i=1}^{n_p}  \uu_{(i)} \sigma_i \vv_{(i)}^\text{T}
\end{equation}
is the SVD  of ${\pmb \Upsilon}_{\dd}$ 
, i.e., $\capu:=[\uu_{(1)} \quad \cdots \quad \uu_{(n_s)}]\in\real^{n_s \times n_s}$, $\capv:=[\vv_{(1)} \quad \cdots \quad \vv_{(n_p)}]\in\real^{n_p \times n_p}$
are orthogonal matrices such that $$\capu^{\rm T}{\pmb \Upsilon}_{\dd}\capv=\Sigma:=\text{diag}(\sigma_1,\ldots,\sigma_{n_p}) \in \real^{n_p \times n_p}$$ 
where $\sigma_1 \geq \sigma_2 \geq \cdots \geq \sigma_{n_p} \geq 0$ are the singular values of ${\pmb \Upsilon}_{\dd}$,  $\uu_{(i)}$ are  the  left singular vectors of ${\pmb \Upsilon}_{\dd}$, and  $\vv_{(i)}$ are the  right singular vectors of ${\pmb \Upsilon}_{\dd}$. 
If ${\pmb \Upsilon}_{\dd}$ does not have full rank, i.e., $\text{rank}({\pmb \Upsilon}_{\dd})=n_r<n_p$, then it holds that $\sigma_{i}=0$ for all $i\in\mathcal{N}=\{n_p-n_r,n_p-n_r+1,\ldots,n_p\}$ (note that $|\mathcal{N}|=n_p-n_r$) and the vectors $\vv_{(i)}$, $i\in\mathcal{N}$, span its nullspace. Hence, said set of vectors is relevant to the identifiability problem since they point to directions of possible unidentifiability because there exist $c_i \neq 0$ such that \eqref{eq:RSSequation} admits solution $$
\Delta \btheta = c_1 \vv_{(n_p-n_r)}+ c_2 \vv_{(n-p-n_r)}+\ldots+c_{n_r} \vv_{(n_p)} \neq \zero_{n_p}.$$ 

Thus,  incorporating all of the above, 
we propose Algorithm \ref{algo.LSERC}, given reference parameters $\btheta^*$, reference input $\uu^*$, sampling times $\{t_0,t_1,t_2,\ldots,t_N\}$, and primary probing directions $D=\{\dd_1,\dd_2,\ldots,\dd_l\}$. There are three stages to the procedure (marked by Stages 1, 2, and 3):
\begin{enumerate}
\item[S1.] \textbf{Primary probing stage}: Complete L-SERC tests for the given primary probing directions $D=\{\dd_i\}$. 
\item[S2.] \textbf{Twin probing stage} ($\epsilon_{\rm twin}>0 \Rightarrow$ ON, $\epsilon_{\rm twin}=0 \Rightarrow$ OFF): Perturb the primary probing directions in one component by adding $\epsilon_{\rm twin} \coord_j$ for some $j$ 
\item[S3.] \textbf{Singularity probing stage} ($\epsilon_{\rm sing}>0 \Rightarrow$ ON, $\epsilon_{\rm sing}=0 \Rightarrow$ OFF): For each of primary and twin probing direction that yields a rank-deficient ${\pmb \Upsilon}_{\dd_i}$, complete an SVD and perturb  $\btheta^*$ in the direction of a linear combination of the singular vectors associated with zero singular values. 
\end{enumerate}

\begin{remark}
Some remarks are in order:
\begin{itemize}
\item If $D=\{\pm \coord_i: i=1,\ldots,n_p\}$ is selected ($|D|=2n_p$), then the primary probing stage checks the natural identifiability test in Definition \ref{defn:identifiability.natural}. 
\item The choice of $\coord_j$ in the twin probing stage can be done systematically for all primary probing directions using, e.g., a Latin hypercube sampling procedure, so as to not repeat the perturbation. 
\item 
The singularity probing stage  (which is inspired by the extended algorithm in \cite{van2022extending}) 
 means the entire algorithm is repeated for new perturbed reference parameters in $\Theta_{sing}$ (if $q>0$), and is repeated again and again etc. (if $q>1$), with the idea being that perturbing along these singular vectors \textit{should} result in no change in the rank but it may because of the so-called singularity problem outlined in \cite{van2022sensitivity} --- see the discussion above concerning the example $y=\max(\theta^3,\theta^5)$. Note that the $\pm$ symbol in Line \ref{eq:algo:SVD} appears so the user explores both directions (because of the possibly piecewise nature of the unknown nonsmooth function) and note that $\epsilon_{\rm sing}=10^{-2} \|\btheta^*\|/\|\Delta \btheta\|$ is chosen in \cite{van2022sensitivity}.
\end{itemize}
\end{remark}
   
In each instance in the algorithm in which the rank of an L-SERC matrix is calculated, this means numerically solving \eqref{eq:1} and \eqref{eq:sens.nonsmooth} simultaneously for identifiability (or \eqref{eq:1} and \eqref{eq:sens.obs} for observability), with $\btheta=\btheta^*$, $\uu=\uu^*$, $\dd_i \in \real^{n_p}$, and calculating ${\pmb \Upsilon}_{\dd_i}$ by sampling the solution at the $N+1$ prescribed times $t=t_k$.

\begin{algorithm}[h!]\caption{Practical L-SERC Algorithm.}
\label{algo.LSERC}
\begin{algorithmic}[1]
\Require $\{(\btheta^*,\uu^*,D=\{\dd_i\},\epsilon_{\rm twin}, \epsilon_{\rm sing}, q\}$ 
\Procedure{}{}  
\State Set $D_{\rm sing}=\emptyset$, $\Theta_{\rm sing}=\emptyset$.\label{algo.jump} 
\For{$i=1,2,\ldots,|D|$}
    \State Calculate $n_i = \text{rank}({\pmb \Upsilon}_{\dd_i})$. \Comment{S1}
    \If {$\epsilon_{\rm sing}>0$ and $n_i<n_p$}
        \State Set $D_{\rm sing} \gets D_{\rm sing} \cup \{\dd_i\}$.
    \EndIf
    \If {$\epsilon_{\rm twin}>0$} \Comment{S2}
        \State Set $\widetilde{\dd}_i \gets \dd_i + \epsilon_{\rm twin} \coord_j$ s.t.\ $\dd_i \neq \coord_j$.  
        \State Calculate $\widetilde{n}_i = \text{rank}({\pmb \Upsilon}_{\widetilde{\dd}_i})$.
        \If {$\epsilon_{\rm sing}>0$ and $\widetilde{n}_i<n_p$}
            \State Set $D_{\rm sing} \gets D_{\rm sing} \cup \{\widetilde{\dd}_i\}$. 
        \EndIf
    \EndIf
\EndFor
\If {$\epsilon_{\rm sing}>0$ and $q>0$} \Comment{S3}
    \ForAll{$\dd_{j} \in D_{\rm sing}$}
        \State Calculate SVD of ${\pmb \Upsilon}_{\dd_j}$.
            \State Set $\btheta^*_{j} \gets \btheta^* \pm \epsilon_{\rm sing} \sum_{k: \sigma_k=0} \vv_{(k)}$.\label{eq:algo:SVD}
            \State Set $\Theta_{\rm sing} \gets \Theta_{\rm sing} \cup \{\btheta^*_j\}$.
    \EndFor
    \State $q \gets q-1$.
    \State Repeat for some $\btheta^* \in \Theta_{\rm sing}$ (i.e., \Goto{algo.jump}).
\EndIf 
\State\Return $\{{\pmb \Upsilon}_{\dd_i}\}$ and $n_i$ for every L-SERC test. 
\EndProcedure
\end{algorithmic}
\end{algorithm} 
  
\begin{example}
We applied Algorithm \ref{algo.LSERC} to the  problem in Example \ref{ex:riot}, without access to the closed-form solution, with $\btheta^*=(1,1)$, $u^*(t)=0$, $D=\{\dd_i=\pm \coord_i: i=1,2\}$, $\epsilon_{\rm twin}=0.01$, $\epsilon_{\rm sing}=0.01$, $q=1$. We obtained the following results (recall $n_p=2$): 
\begin{align*}
&\mathrm{S1}: \dd_i = -\coord_1, \coord_2 \Rightarrow  \mathrm{rank}({\pmb \Upsilon}_{\dd_i})=2=n_p.\\
&\mathrm{S1}: \dd_i = \coord_1,-\coord_2 \Rightarrow  \mathrm{rank}({\pmb \Upsilon}_{\dd_i})=1<n_p.\\
&\mathrm{S2}: \widetilde{\dd}_i = -\coord_1+\epsilon_{\rm twin}\coord_2, \coord_2+\epsilon_{\rm twin}\coord_1 \Rightarrow  \mathrm{rank}({\pmb \Upsilon}_{\widetilde{\dd}_i})=2.\\
&\mathrm{S2}: \widetilde{\dd}_i = \coord_1+\epsilon_{\rm twin}\coord_2,-\coord_2+\epsilon_{\rm twin}\coord_1 \Rightarrow  \mathrm{rank}({\pmb \Upsilon}_{\widetilde{\dd}_i})=1.\\
&\mathrm{S3}: \mathrm{SVD}({\pmb \Upsilon}_{\dd_i \in D_{\rm sing}}) \Rightarrow  \sigma_2=0, \vv_{(2)}=\begin{bsmallmatrix} -1 \\ 0 \end{bsmallmatrix}.
\\
&\qquad \Rightarrow \btheta^*=\begin{bsmallmatrix} 1\pm \epsilon_{\rm sing} \\ 1 \end{bsmallmatrix}=\begin{bsmallmatrix} 1.01 \\ 1 \end{bsmallmatrix}, \begin{bsmallmatrix} 0.99 \\ 1 \end{bsmallmatrix}.\\
&\mathrm{S1} (\btheta^*=\begin{bsmallmatrix} 1.01 \\ 1 \end{bsmallmatrix}): \dd_i = \pm \coord_1, \pm \coord_2 \Rightarrow  \mathrm{rank}({\pmb \Upsilon}_{\dd_i})=1.\\
&\mathrm{S1} (\btheta^*=\begin{bsmallmatrix} 0.99 \\ 1 \end{bsmallmatrix}): \dd_i = \pm \coord_1, \pm \coord_2 \Rightarrow  \mathrm{rank}({\pmb \Upsilon}_{\dd_i})=2.
\end{align*}  
 Hence, \eqref{riot1} is partially identifiability but not naturally identifiable. Based on the findings, we hypothesize (i) identifiability in the third quadrant (relative to $\btheta^*=(1,1)$); (ii) unidentifiability in the second quadrant; and (iii) $\btheta^*=(1,1)$ lies on a nonsmoothness. And we would be correct; we know from the closed-form solution that the ``region of identifiability'' is given by $\Theta_{\rm ident}=\{\btheta:\theta_1 \leq \theta_2\}$, while  $\Theta_{\rm unident}=\{\btheta:\theta_1 > \theta_2\}$.
\end{example}

\section{Application: Stommel-Box Climate Model}\label{sec:stommel} 
The ocean's role in the climate system is significant as it transports and stores a substantial amount of heat. However, oceans have been warming at an alarming rate, which can have crucial consequences on the ecosystem and life in general on the planet Earth; the North Atlantic is where this crisis is manifesting itself the most -- see \cite{li2021observation} and references therein. Thus, many studies and models -- such as the Stommel-box model \cite{stommel1961thermohaline} -- focus on studying the relation between heat flux, freshwater flux, and their direct effect on the ocean circulation to understand and predict causes and consequences of climate change, due to warming (increased heat input) and changes in freshwater/salinity concentration/distribution due to, for example, the melting of ice sheets and permafrost layers in the Arctic \cite{vincent2020arctic}. Said Stommel-box model consists of two main boxes (interactions): a high altitude one (e.g., the surface of the ocean) and a low altitude one (e.g., the bottom of the ocean). Normally, near the surface of the ocean, there is a higher temperature and lower salinity (this is reversed near the bottom of the ocean). In \cite{budd2022dynamic}, the authors studied a nonsmooth version of the Stommel-box model, and demonstrated that it was able to predict tipping point occurrences better than the smooth versions available in the literature. Said tipping points are crucial, since they mark dramatic and irreversible changes taking place, indicating the start of a new phase of the dynamics (which could be related to, e.g., the beginning of a disastrous climate event). The aforementioned nonsmooth Stommel-box model is given as follows:
\begin{subequations}\label{eq:Stommel_model}
\begin{align}
            \dot{T}&=\theta_1 + u_1 -T-T|T-V|,& &T(0)=T_0,  \\
        \dot{V}&=\theta_2 + u_2 -V\theta_3-V|T-V|,& &V(0)=V_0,
\end{align}
\end{subequations}
where $T$ and $V$ represent the dimensionless temperature and salinity, $t$ is the dimensionless time, the parameters $\theta_1$, $\theta_2$ and $\theta_3$ represent the strength of the thermal forcing, the strength of the freshwater forcing, and the ratio of the thermal and freshwater surface restoring time scales, respectively. The inputs $u_1$ and $u_2$ represent variations in parameters from their baseline values $\theta_1$ and $\theta_2$;
following \cite{budd2022dynamic}, we assume that 
$u_1(t)=B \sin(\Omega t)$ and $u_2(t)=\hat{B} \sin(\Omega t)$, where $B$, $\hat{B}$, and $\Omega$ are constants. The model in \eqref{eq:Stommel_model} is relatively simple but is able to capture important features of the dynamics and provide important information; for instance, the bifurcation diagram obtained from this model matches the bifurcation studies for two-dimensional and three-dimensional ocean models \cite{dijkstra2013nonlinear}.

\begin{figure} [!htb]
	\centering
 \includegraphics[width=1\linewidth]{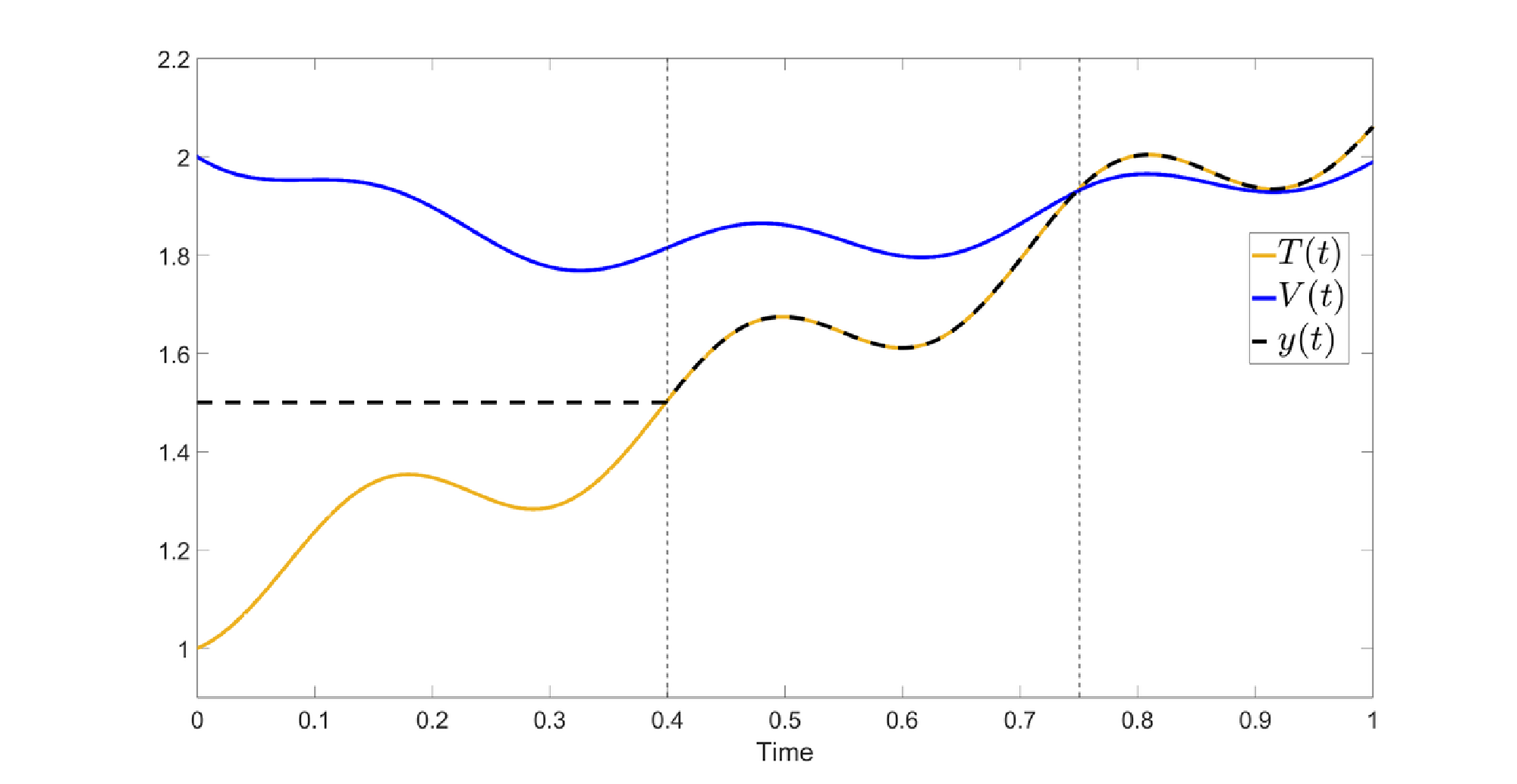}
	\caption{Trajectories of the state variables $T$, $V$, and the output $y$. 
    The output $y$ crosses a nonsmooth threshold when $T(t)=T_{\min}$ (indicated by the first vertical dashed line). 
	The second vertical dashed line corresponds to a nonsmooth threshold in the model being crossed  when $T(t)=V(t)$.}
 \label{fig:Stommel Solution vs Output}
\end{figure}

In order to apply the model in \eqref{eq:Stommel_model} in an effective way for different kinds of studies, one may need to determine the parameters $\theta_1$, $\theta_2$ and $\theta_3$ that best fit certain environments; for that, one will need to take measurements of available information (e.g., measurement of temperature and/or salinity) to identify the model's parameters. However, identifiability and observability methods available in the literature cannot be used for the model in \eqref{eq:Stommel_model} due to the nonsmoothness present in the RHS. Hence, we use the newly developed L-SERC methods in Section \ref{sec:nonsmoothSERC} to analyze identifiability and observability of the nonsmooth Stommel-box model in \eqref{eq:Stommel_model}. From many possible output functions, we assume we have access to temperature measurements. Hence, we take $y=\max(T,T_{\min})$, where $T_{\min}$ is the lowest measurable temperature due to sensory, equipment, and/or computational limitations (see \cite{tyler2013using} for example). 
We consider initial conditions $T_{0}=1,V_{0}=2$, reference control inputs $\uu^*(t)=(u_1^*(t),u_2^*(t))=(2\sin(20t),\sin(20t))$ (i.e., $B=2$, $\hat{B}=1$, $\Omega=20$), reference parameters $\btheta^*=(\theta_1^*,\theta_2^*,\theta_3^*)=(3,1.1,0.3)$, and lowest measurable temperature as $T_{\min}=0.5$. The values of the parameters and initial conditions considered are within the range of values used in relevant literature such as \cite{budd2022dynamic}. 
Figure \ref{fig:Stommel Solution vs Output} depicts the trajectories of the states and the output. 

\begin{figure} [!htb]
    \centering
\includegraphics[width=1\linewidth]{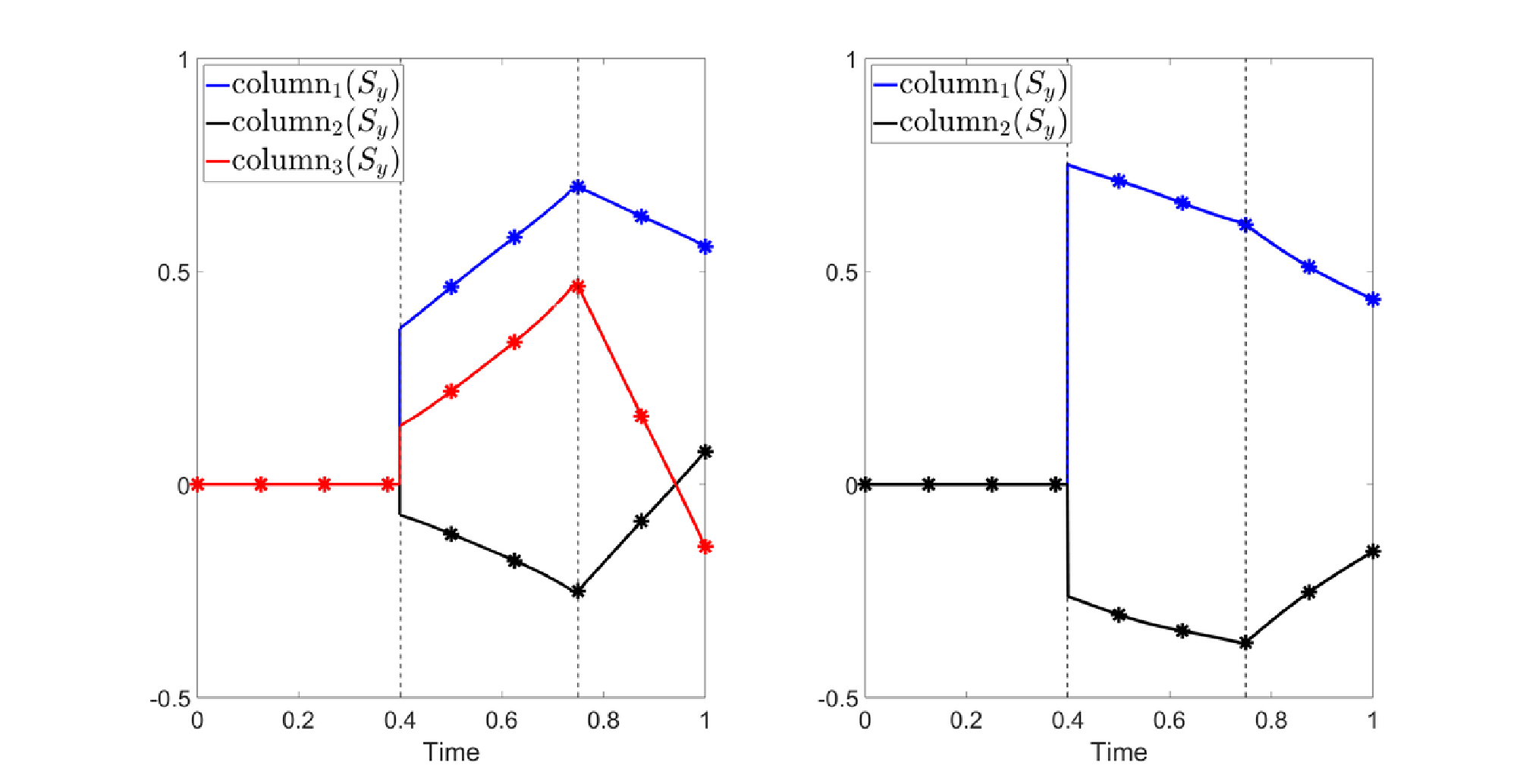}
 
 \caption{This figure presents the output sensitivity function $\caps_{\y}^{\rm L}(t)=\mathrm{lshift}(\capy^*(t)) \in \real^{1 \times 3}$, with $\capy^*(t)$ from \eqref{eq:Stommel_X_dot} for identifiability (left panel), and $\caps_{\y}^{\rm L}(t)=\mathrm{lshift}(\capy^*(t)) \in \real^{1 \times 2}$, with $\capy^*(t)$ from \eqref{eq:Stommel_X_dot_obs} for observability (right panel). 
	When computing $\caps_{\y}(t)$ for both the identifiability and observability, we can see that, as shown in Figure \ref{fig:Stommel Solution vs Output}, the nonsmoothness happens in the output when $T(t)=T_{\min}$ (first dashed line) and in the model when $T(t)=V(t)$ (second dashed line). The stars in this figure represent the sample times $\{t_k\}$ used to construct ${\pmb \Upsilon}_{\dd}$ in \eqref{LSERC_ID}.} 
 \label{fig:Stommel S_Y}
\end{figure}

First we tested identifiability: 
in order for us to   apply Algorithm 1, we first obtain the nonsmooth sensitivity system in \eqref{eq:sens.nonsmooth} for the Stommel-box model in \eqref{eq:Stommel_model}, which can be derived using the LD-derivative rules as 
\begin{equation}\label{eq:Stommel_X_dot}
\begin{aligned} 
&\dot{\capx}_1 =[d_1 \quad 1 \quad 0 \quad 0]
    -(1+|T^*-V^*|)\capx_1\\
    &\qquad     -T^*\mathrm{fsign}(T^*-V^*,\capx_1-\capx_2) \; (\capx_1-\capx_2),\\
&\dot{\capx}_2 =[d_2 \quad 0 \quad 1 \quad 0]- V^* [d_3 \quad 0 \quad 0 \quad 1]\\
    &\qquad -(\theta_3^*+|T^*-V^*|)\capx_2\\
    &\qquad  -V^*\mathrm{fsign}(T^*-V^*,\capx_1-\capx_2) \; (\capx_1-\capx_2), \\
&\capy = {\rm \bf slmax}
([T^* \quad \capx_1],[T_{\min} \quad \bm{0}_{1 \times 4}])
,\\
&\capx(0)=\bm{0}_{2 \times 4},
\end{aligned}
\end{equation}
where $\capx_1(t):=\mathrm{row}_1(\capx(t)) \in \real^{1 \times 4}$ and $\capx_2(t):=\mathrm{row}_2(\capx(t)) \in \real^{1 \times 4}$. We applied Algorithm 1 in the primary probing stage to assess local identifiability in the natural sense (i.e., taking $D=\{d_i=\pm \coord_i: i=1,2,3\}$). The left panel of Figure \ref{fig:Stommel S_Y} depicts the output sensitivity function $\caps_{\y}$ w.r.t.\ every parameter, along with the samples (in stars) taken to construct the L-SERC matrix 
$${\pmb \Upsilon}_{\dd}
=\mathrm{lshift}
(\capy^*(t_0),\ldots,\capy^*(t_9))
\in \real^{10 \times 3},$$
with sampling times $\{t_0,t_1,t_2,t_3,t_4,t_5,t_6,t_7,t_8,t_9\}=\{0,0.125,0.25,0.375,0.5,0.625,0.75,0.875,1\}$. We obtained full rank L-SERC matrices in each $\dd_i=\pm \coord_i$ case, which implies natural identifiability. We then executed the twin probing stage with $\epsilon_{twin}=0.01$ in Algorithm 1 to further explore this problem and solidify our findings; again the L-SERC matrices we computed had full rank, from which we hypothesize full structural identifiability in this model. Because $D_{\rm sing}=\emptyset$, the singularity probing stage was skipped.

We also tested the model for observability. For that,  the parameters are the initial conditions of the states, i.e., $\f_0(\btheta)=\btheta=(T_0,V_0)$, and the relevant nonsmooth sensitivity system in \eqref{eq:sens.obs} takes the following form:
\begin{equation}\label{eq:Stommel_X_dot_obs}
\begin{aligned} 
&\dot{\capx}_1 =    -(1+|T^*-V^*|)\capx_1\\ 
&\qquad -T^*\mathrm{fsign}(T^*-V^*,\capx_1-\capx_2) \; (\capx_1-\capx_2),\\
&\dot{\capx}_2 =  
    -(\theta_3^*+|T^*-V^*|)\capx_2\\
    &\qquad -V^*\mathrm{fsign}(T^*-V^*,\capx_1-\capx_2) \; (\capx_1-\capx_2), \\
&\capy = {\rm \bf slmax}
([T^* \quad \capx_1],[T_{\min} \quad \bm{0}_{1 \times 4}]),\\
&\capx(0)=[\dd \quad \II_3],
\end{aligned}
\end{equation}
The right panel of  Figure \ref{fig:Stommel S_Y} depicts the output sensitivity functions w.r.t.\ the two initial conditions; we applied the primary probing stage in Algorithm 1, i.e., $D=\{\dd_i=\pm \coord_i: i=1,2\}$, and the twin probing stage (with $\epsilon_{twin}=0.01$ ) and, based on the samples, we conclude natural observability of the system. 


We note here that the introduced approach is shown to be effective in assessing identifiability and observability in two significant ways.  First, this L-SERC approach is applicable without any \textit{a priori} knowledge of the input-output system encountering nonsmoothness (or needing to check for it along the way); this can be seen clearly in the simulation results handling nonsmoothness in the model and the output.
Second, the new L-SERC test is capable of providing information on identifiability and observability even if during the simulation/experiment there are intervals of time in which the output is not sensitive to the parameters; this can be seen in Figure \ref{fig:Stommel S_Y}, where the output does not provide information for identifiability and/or observability for $t\in[0,0.4]$, however, sampling the L-SERC matrix beyond the first nonsmoothness allowed for the assessment of identifiability and/or observability since enough samples were collected after $t=0.4$.


\section{Concluding Remarks}\label{sec:conclusion}
This paper provides a novel theory, which is practically implementable, for testing local identifiability and observability of \textit{nonsmooth} dynamical and control systems. The newly developed theory is based on  lexicographic directional   derivatives, which are used to create a first-order Taylor-like approximation theory for a broad class of nonsmooth functions. Said ``nonsmooth linearization'' approximation theory enabled us to derive a nonsmooth lexicographic-based sensitivity rank condition (L-SERC) test that is analogous to, and recovers, the smooth SERC test. 
We introduced new definitions of identifiability/observability that fit the nonsmooth analysis environment, and put forward an algorithm that incorporates many advancements that took place in the smooth SERC literature. 
We applied the new algorithm to a recent nonsmooth climate model to demonstrate the effectiveness of the proposed methods; the numerical application we provided here makes a strong case for how our approach is valuable in assessing identifiability and observability even near/on/about nonsmooth points, without the need to check them or conduct any smoothing approximations (which can only make climate models less accurate as observed in \cite{budd2022dynamic}).

The results of this paper, we hope, will encourage research efforts aimed at analyzing nonsmooth or discontinuous dynamical systems   properties  in a direct way without the need for smoothing approximations and its respective inaccuracies. 
A nonsmooth sensitivity theory has been established for other nonsmooth dynamical systems using lexicographic directional derivatives, such as differential-algebraic equations \cite{2016_Stechlinski_Barton_LD_DAEs,stechlinski2021nonsmooth} and optimization-constrained ODEs \cite{2018_DAEOs}, among others. Thus, it should be possible to extend the L-SERC identifiability and observability theory in the present paper to said systems. 
Another direction for future work is to extend the sensitivity-analysis based controllability test (see, e.g., \cite{van2022sensitivity}) to the nonsmooth setting by establishing adjoint sensitivity theory for nonsmooth ODEs. 
Lastly, we note that a ``negative L-SERC test'' is desirable (i.e., one that determines when system parameters are not identifiable), but this is challenging even in the smooth setting (see again the discussion about isolated output sensitivity singularities in Section \ref{sec:practicalalgo}). In the nonsmooth setting, things prove yet more complicated because of available necessary and sufficient nonsmooth inverse function theory (see, e.g., \cite[Theorem 3]{gowda2004inverse}) and because the L-subdifferential may be a strict superset of  Clarke's generalized derivative in the general L-smooth setting.

  \section*{Author Contributions}
The order of authors' names reflects their overall contributions to the paper: {\bf Peter Stechlinski}: Conceptualization, Analysis, Writing (original draft, review \& editing), Supervision. {\bf Sameh Eisa}: Conceptualization, Analysis, Writing (original draft, review \& editing), Supervision. {\bf Hesham Abdelfattah}: Literature Review, Simulations, Writing (original draft).

\appendices

  \section{Proof of Theorem \ref{thm.firstorder.LDderivative}}\label{sec.proof1} 

Before proving Theorem \ref{thm.firstorder.LDderivative}, we need a relevant lemma.

\begin{lemma}\label{lemma.homogenization}
 Let $X \subseteq \real^n$ be open and $\f:X \to \real^n$ be L-smooth at $\xnot \in X$. Then, given any $\x \in \real^n$,   
\begin{equation}\label{eq:homogenizationrelation}
[\f^{(1)}(\coord_1)  \quad \cdots \quad \f^{(n)}(\coord_n)] \; (\x-\xnot)=\f'(\xnot;\x-\xnot),
\end{equation}
 where $\coord_i$ is the $i^{\rm th}$ standard basis vector in $\real^n$. 
\end{lemma}
 
 \begin{pf}
 Let $\M=[\Delta \x \quad \II_n] \in \real^{n \times (1+n)}$ where $\Delta \x=\x-\xnot$. Then, by definition of the LD-derivative in \eqref{eq:LDderivative},
\begin{align}
&\f'(\xnot;[\Delta \x \quad \II_n]) \begin{bsmallmatrix} 0 \\ \Delta \x \end{bsmallmatrix}\notag\\
&=[\f^{(0)}(\Delta \x) \quad \f^{(1)}(\coord_1) \quad \cdots \quad \f^{(n)}(\coord_n)]\begin{bsmallmatrix} 0 \\ \Delta \x \end{bsmallmatrix}\notag\\
&=[\f^{(1)}(\coord_1) \quad \cdots \quad \f^{(n)}(\coord_n)]\Delta \x.\label{eq:homoglemma:1}
\end{align}
On the other hand, since $\M$ has full row rank,  \eqref{eq:LDderivative:vs:Lderivative} gives  
$$
\f'(\xnot;[\Delta \x \quad \II_n])
=\JJ_{\rm L}\f(\xnot;[\Delta \x \quad \II_n])[\Delta \x \quad \II_n]$$
so that
\begin{align}
&\f'(\xnot;[\Delta \x \quad \II_n]) \begin{bsmallmatrix} 0 \\ \Delta \x \end{bsmallmatrix}\notag\\
&=\JJ_{\rm L}\f(\xnot;[\Delta \x \quad \II_n])[\Delta \x \quad \II_n]\begin{bsmallmatrix} 0 \\ \Delta \x \end{bsmallmatrix}\notag\\
&=\JJ_{\rm L}\f(\xnot;[\Delta \x \quad \II_n])\Delta \x.\label{eq:homoglemma:2}
\end{align}
Hence, from \eqref{eq:homoglemma:1} and \eqref{eq:homoglemma:2}, 
$$[\f^{(1)}(\coord_1) \quad \cdots \quad \f^{(n)}(\coord_n)]\Delta \x=\JJ_{\rm L}\f(\xnot;[\Delta \x \quad \II_n])\Delta \x.$$
Then, by choosing $\x=\xnot+\coord_1$, $\x=\xnot+\coord_2$, etc., so that $\Delta \x=\coord_1$, $\Delta \x=\coord_2$, etc., it follows that 
\begin{equation}\label{eq:homogenizatino:conc1}
[\f^{(1)}(\coord_1) \quad \cdots \quad \f^{(n)}(\coord_n)]=\JJ_{\rm L}\f(\xnot;[\Delta \x \quad \II_n]).
\end{equation}
By definition of the L-derivative in Equation \eqref{eq:Lderivative}, $$\JJ_{\rm L}\f(\xnot;[\Delta \x \quad \II_n])=\JJ \f^{(n)}(\zero_n),$$ 
with $\f^{(n)}$ being linear since $[\Delta \x \quad \II_n]$ spans $\real^n$, i.e.,  $\f^{(n)}:\x \mapsto \mathbf{A} \x$
with $\mathbf{A}=\JJ_{\rm L}\f(\xnot;[\Delta \x \quad \II_n])$. 
On the other hand, recalling that $\f^{(0)}(\Delta \x)=\f'(\xnot;\Delta \x)$ by definition of the homogenization sequence in \eqref{eq:Lsmoothfunctions},   
\begin{equation}\label{eq:homogenizatino:conc3}
\f'(\xnot;\Delta \x)=\f^{(0)}(\Delta \x)=\f^{(n)}(\Delta \x)
\end{equation}
where the final equality comes from Lemma 2.1 Property 4 in \cite{khan2014generalized}. 
It therefore follows from \eqref{eq:homogenizatino:conc1}--\eqref{eq:homogenizatino:conc3} that  
\begin{align*}
\f'(\xnot;\Delta \x)
&=\f^{(n)}(\Delta \x)=\JJ_{\rm L}\f(\xnot;[\Delta \x \quad \II_n]) \Delta \x \\
&=[\f^{(1)}(\coord_1) \quad \f^{(2)}(\coord_2) \quad \cdots \quad \f^{(n)}(\coord_n)]\Delta \x. 
\end{align*} 
and hence the result follows. \qed
 \end{pf} 
    


\begin{pf}\textbf{[Proof of Theorem \ref{thm.firstorder.LDderivative}]} Following the proof of \cite[Theorem 3.1.2]{Scholtes},  suppose  not, i.e., 
$$\| \f(\xnot+\dd_n)-\f(\xnot)
-\f'(\xnot;[\dd_n \quad \II_n]) \begin{bsmallmatrix} 0 \\ \dd_n \end{bsmallmatrix}   \|/\|\dd_n\| \not\to 0$$
for some sequence $\dd_n \to \zero_n$. Assume, without loss of generality, that $\dd_n/\|\dd_n\|\to \bar{\dd}\in \{\x \in \real^n:\|\x\|=1\}$ (since, if not, we can choose a suitable subsequence). 
Letting $\alpha_n=\|\dd_n\|$, observe that 
\begin{align}
&0 \leq \| \f(\xnot+\dd_n)-\f(\xnot)-\f'(\xnot;[\dd_n \quad \II_n]) \begin{bsmallmatrix} 0 \\ \dd_n \end{bsmallmatrix}  \| / \|\dd_n\|\notag\\
&=
\tfrac{1}{\alpha_n} 
\|
\f(\xnot+\dd_n)
-\f(\xnot+\alpha_n \bar{\dd})\notag\\
&+\f'(\xnot;[\alpha_n \bar{\dd} \quad \II_n])\begin{bsmallmatrix} 0 \\ \alpha_n \bar{\dd} \end{bsmallmatrix}
-\f'(\xnot;[\dd_n \quad \II_n])\begin{bsmallmatrix} 0 \\ \dd_n \end{bsmallmatrix}\notag\\
&+\f(\xnot+\alpha_n \bar{\dd})
-\f(\xnot)
-\f'(\xnot;[\alpha_n \bar{\dd} \quad \II_n])\begin{bsmallmatrix} 0 \\ \alpha_n \bar{\dd} \end{bsmallmatrix}  
\|\notag\\
&\leq
\tfrac{1}{\alpha_n}
\|\f(\xnot+\dd_n)-\f(\xnot+\alpha_n \bar{\dd})\|\label{eq:expression1}\\
&+
\tfrac{1}{\alpha_n}
\|\f'(\xnot;[\alpha_n \bar{\dd} \quad \II_n])\begin{bsmallmatrix} 0 \\ \alpha_n \bar{\dd} \end{bsmallmatrix}
-\f'(\xnot;[\dd_n \quad \II_n])\begin{bsmallmatrix} 0 \\ \dd_n \end{bsmallmatrix}\|\label{eq:expression2}\\
&+
\tfrac{1}{\alpha_n}
\|\f(\xnot+\alpha_n \bar{\dd})-\f(\xnot)-\f'(\xnot;[\alpha_n \bar{\dd} \quad \II_n])\begin{bsmallmatrix} 0 \\ \alpha_n \bar{\dd} \end{bsmallmatrix}\|.\label{eq:expression3}
\end{align} 
We proceed by ``squeezing'' the expressions in \eqref{eq:expression1}--\eqref{eq:expression3}. For \eqref{eq:expression1} and \eqref{eq:expression2}, we follow the same argument as in \cite[Theorem 3.1.2]{Scholtes}: First, by local Lipschitz continuity of $\f$, and since $\alpha_n \to 0$, there exist $N_1 \in \posint$ and a constant $L$ such that, for all $n \geq N_1$,
\begin{align*}
\|\f(\xnot+\dd_n)-\f(\xnot+\alpha_n \bar{\dd})\| 
& \leq L\|\dd_n-\alpha_n \bar{\dd}\|.
\end{align*}
Next, in considering \eqref{eq:expression2},  we note that
\begin{align*}
&\f'(\xnot;[\alpha_n \bar{\dd} \quad \II_n])\begin{bsmallmatrix} 0 \\ \alpha_n \bar{\dd}-\dd_n \end{bsmallmatrix}\\
&=
[\f^{(0)}(\alpha_n \bar{\dd}) \quad \f^{(1)}(\coord_1) \quad \cdots \quad \f^{(n)}(\coord_n)] \begin{bsmallmatrix} 0 \\ \alpha_n \bar{\dd}-\dd_n \end{bsmallmatrix}\\
&=[\f^{(1)}(\coord_1) \quad \cdots \quad \f^{(n)}(\coord_n)] (\alpha_n \bar{\dd}-\dd_n),
\end{align*}
where $\coord_i$ is the $i^{\rm th}$ standard basis vector in $\real^n$.   
Hence, 
\begin{align*}
&\|\f'(\xnot;[\alpha_n \bar{\dd} \quad \II_n])\begin{bsmallmatrix} 0 \\ \alpha_n \bar{\dd} \end{bsmallmatrix}-\f'(\xnot;[\dd_n \quad \II_n]) \begin{bsmallmatrix} 0 \\ \dd_n \end{bsmallmatrix}\|\\
&=
\| [\f^{(1)}(\coord_1) \quad \cdots \quad \f^{(n)}(\coord_n)] (\alpha_n \bar{\dd}-\dd_n)\|\\
&\leq 
\| [\f^{(1)}(\coord_1) \quad \cdots \quad \f^{(n)}(\coord_n)]\| \| \alpha_n \bar{\dd}-\dd_n\|.
\end{align*}
Lastly, in considering \eqref{eq:expression3}, observe that 
\begin{align*}
&\|\f(\xnot+\alpha_n \bar{\dd})-\f(\xnot)-\f'(\xnot;[\alpha_n \bar{\dd} \quad \II_n]) \begin{bsmallmatrix} 0 \\ \alpha_n \bar{\dd} \end{bsmallmatrix}\| \\
&=
\|\f(\xnot+\alpha_n \bar{\dd})-\f(\xnot)-[\f^{(1)}(\coord_1) \; \cdots \; \f^{(n)}(\coord_n)] \alpha_n \bar{\dd}\|. 
\end{align*}

Collecting the inequalities, we have that, for $n \geq N_1$,
\begin{align*}
&\| \f(\xnot+\dd_n)-\f(\xnot)-\f'(\xnot;[\dd_n \quad \II_n]) \begin{bsmallmatrix} 0 \\ \dd_n \end{bsmallmatrix} \| / |\dd_n\| \\
&\leq
\tfrac{1}{\alpha_n}
(
L\|\dd_n-\alpha_n \bar{\dd}\|\\
& +
\| [\f^{(1)}(\coord_1) \; \cdots \; \f^{(n)}(\coord_n)]\| \| \alpha_n \bar{\dd}-\dd_n\|\\
& +
\|\f(\xnot+\alpha_n \bar{\dd})-\f(\xnot)-[\f^{(1)}(\coord_1) \; \cdots \; \f^{(n)}(\coord_n)] \alpha_n \bar{\dd}\|)\\
&=
(L+\| [\f^{(1)}(\coord_1) \; \cdots \; \f^{(n)}(\coord_n)]\|) \| \tfrac{\dd_n}{\alpha_n}-\bar{\dd}\|\\
&\; +
\left\|\tfrac{\f(\xnot+\alpha_n \bar{\dd})-\f(\xnot)}{\alpha_n}-[\f^{(1)}(\coord_1) \;  \cdots \; \f^{(n)}(\coord_n)] \bar{\dd} \right\|.
\end{align*} 
But
$(L+\| [\f^{(1)}(\coord_1) \; \cdots \; \f^{(n)}(\coord_n)]\|)  \| \tfrac{\dd_n}{\alpha_n}-\bar{\dd} \| \to 0$
as $n\to\infty$ since $\dd_n/\alpha_n \to \bar{\dd}$, and 
\begin{align*} 
\|\tfrac{\f(\xnot+\alpha_n \bar{\dd})-\f(\xnot)}{\alpha_n}-[\f^{(1)}(\coord_1) \; \cdots \; \f^{(n)}(\coord_n)] \bar{\dd} \|
&\to 0
\end{align*} 
as $n \to \infty$ by Lemma \ref{lemma.homogenization} with $\x=\bar{\dd}+\xnot$. Hence, we get 
$$\lim_{\dd_n \to \zero_n} \tfrac{\| \f(\xnot+\dd_n)-\f(\xnot)-\f'(\xnot;[\dd_n \quad \II_n])\begin{bsmallmatrix} 0 \\ \dd_n \end{bsmallmatrix}\|}{\|\dd_n\|} = 0,$$ 
which is a contradiction. \qed
\end{pf}

\section*{References}

\bibliographystyle{IEEEtranS}

\bibliography{LDSERC}

\end{document}